\documentclass[final, 3p, 12pt]{elsarticle}

\usepackage{algorithm}
\usepackage{algorithmicx}
\usepackage{algpseudocode}
\makeatletter \renewcommand*{\ALG@name}{Method}

\usepackage{amssymb}
  \usepackage{amsmath,eqnarray}
  \usepackage{amsthm}
   
   \theoremstyle{plain}
  \usepackage{bm}
  
\newcommand\+{\mkern2mu}
  \usepackage{verbatim}
  \usepackage{textcomp} 
 \usepackage{longtable}
\usepackage{booktabs} 
\usepackage{dcolumn}  
\usepackage{multirow} 
\usepackage{lettrine}
\usepackage{simplemargins}
 \setallmargins{.8in}
\newtheorem{thm}{Theorem}
\makeatletter
\renewcommand*\env@matrix[1][c]{\hskip -\arraycolsep
  \let\@ifnextchar\new@ifnextchar
  \array{*\c@MaxMatrixCols #1}}
\makeatother

\newtheorem{cor}{Corollary}
\newtheorem{lem}{Lemma}

\newtheorem{rem}{Remark}
\theoremstyle{definition}
\newtheorem{defn}{Definition}
\newtheorem*{pf}{Proof}
\newcommand{\oant}{\mbox{OA}(N,k,2,t)}
\newcommand{\oan}{\mbox{OA}(N,k,s,t)}
\newcommand{\ca}{\mbox{CA}^*_{\lambda}(k,s,t)}
\newcommand{\pa}{\mbox{PA}^*_{\lambda}(k,s,t)}

 \newcommand{\Giso}{G^{\rm iso}(k,2)}
\newcommand{\M}{{\bf M}}
\newcommand{\PP}{{\bf P}}

\newcommand{\ppi}{{\boldsymbol \Pi}}
\newcommand{\I}{{\bf I}}
\newcommand{\1}{{\bf 1}}
\newcommand{\0}{{\bf 0}}

\newcommand{\bbb}{{\bf{b}}}
\newcommand{\bb}{{\bf{b}}}
\newcommand{\cc}{{\bf c}}
\newcommand{\uu}{{\bf u}}
\newcommand{\uuuu}{\hat {{\bf u}}}

\newcommand{\x}{{\bf x}}
\newcommand{\vv}{{\bf v}}
\newcommand{\vvvv}{{\hat{\bf v}}}
\newcommand{\V}{{\bf V}}
\newcommand{\y}{{\bf y}}
\newcommand{\s}{{\bf s}}
\newcommand{\e}{{\bf e}}
\newcommand{\z}{{\bf z}}

\newcommand{\dd}{{\bf d}}

\newcommand{\Y}{{\bf Y}}
\newcommand{\R}{{\bf R}}
\newcommand{\U}{{\bf U}}
\newcommand{\ZZ}{{\bf Z}}
\newcommand{\E}{{\bf E}}
\newcommand{\X}{{\bf X}}
\newcommand{\w}{{\bf w}}

\newcommand{\B}{{\bf B}}
\newcommand{\Glp}{G^{\rm LP}}
\newcommand{\Glpp}{G(\A,\bb,\B, \dd,\cc)^{\rm LP}}
\newcommand{\Gabc}{G(\B,\ddd, \cc)}
\newcommand{\ddd}{{\bf d}}
\newcommand{\D}{{\bf D}}

\newcommand{\A}{{{\bf A }}}

\newcommand{\uuu}{{\bf u}}
\newcommand{\GABocd}{G^{\rm Null}_{(\A,\B,\ddd,\cc)}}

\newcommand{\GABbcdn}{G(\A,\bb,\B, \dd,\cc)}
\usepackage{url}
\begin{document}
\begin{frontmatter}
\title{Finding the symmetry group of an LP with equality constraints and its application to classifying orthogonal arrays}
\author[AFIT1]{Andrew J.~Geyer}
\ead{andrew.geyer@afit.edu}
\author[AFIT1]{Dursun A.~Bulutoglu\corref{cor1}}
\ead{dursun.bulutoglu@gmail.com}
\author[WVU]{Kenneth J.~Ryan}
\ead{kjryan@mail.wvu.edu}
\address[AFIT1]{Department of Mathematics and Statistics, Air Force Institute of Technology,\\Wright-Patterson Air Force Base, Ohio 45433, USA}
\address[WVU]{Department of Statistics, West Virginia University,\\ Morgantown, West Virginia 26506, USA}
\cortext[cor1]{Corresponding author}
\journal{Discrete Optimization}
\begin{abstract}
For a given  linear program (LP) a permutation of 
its variables that sends feasible points to
feasible points and preserves the objective function value of each of its feasible points is a symmetry  of the LP. 
The set of all symmetries of an LP, denoted by $\Glp$, is the symmetry group of the LP.
Margot [F. Margot, 50 Years of Integer Programming 1958-2008 (2010), 647-686] described a method for computing  a subgroup of the symmetry group $\Glp$  of an LP.
This method computes $\Glp$ when the LP has only non-redundant inequalities and its feasible set satisfies no equality constraints.  However, when the feasible set of the LP satisfies equality constraints this method
finds only a subgroup of $\Glp$ and can miss symmetries. We  develop a method for finding the symmetry group of a feasible LP whose feasible set satisfies equality  constraints. We  apply  this method to find and exploit the previously unexploited symmetries of an orthogonal array defining integer linear program (ILP)  within the branch-and-bound (B\&B) with isomorphism pruning 
algorithm  [F. Margot, Symmetric ILP: Coloring and small integers, Discrete Optimization 4 (1) (2007), 40-62]. Our method reduced the running time for finding all OD-equivalence  
classes of OA$(160,8,2,4)$ and OA$(176,8,2,4)$ by factors of $1/(2.16)$ and $1/(1.36)$ compared to the fastest known method [D. A. Bulutoglu and K. J. Ryan, Integer programming for classifying orthogonal arrays,  Australasian Journal of Combinatorics 70 (3) (2018), 362-385]. These were the two bottleneck cases that could not have been  solved until the B\&B  with isomorphism pruning algorithm was applied. 
 Another key finding of this paper is that  
converting inequalities to equalities by introducing slack variables and exploiting  the symmetry group of the resulting ILP's LP relaxation within the B\&B with isomorphism pruning algorithm can reduce the computation time by several orders of magnitude when enumerating a set of all non-isomorphic solutions of an ILP.
\end{abstract}
\end{frontmatter}
\section{Introduction}
\label{sec:intro}
A branch-and-bound (B\&B) algorithm can be used to find an optimum solution or enumerate all optimum solutions to an integer linear program (ILP) of the form
 \begin{equation}\label{eqn:geneqILP}
\begin{array}{rl}
\min &   \cc^{\top}\x \\
\mbox{s.t.} &  \A \x=\bb, \quad \x \in \mathbb Z^n,\\
&  \B\x\leq \ddd.
\end{array}
\end{equation}
 Let $\x$ be called a {\em partial solution} of ILP~(\ref{eqn:geneqILP}) if each element of a strict subset of 
 entries of  $\x$ has been assigned integer values and the remaining entries are not fixed.
One way a  B\&B algorithm that 
 branches on the integer values of the variables of $\x$ prunes a partial solution, i.e., a  node of its backtrack search tree, is by infeasibility (pruning by infeasibility). 
 A partial solution is pruned by infeasibility if the linear programming (LP) relaxation of the {\em subproblem} created 
 from ILP~(\ref{eqn:geneqILP}) by assigning the fixed integer values in the partial solution $\x$ to their corresponding variables is proven to be infeasible by solving the LP relaxation. 
Another way of pruning is by comparing the optimum LP relaxation value of a subproblem to that of the best known solution. If this objective function value is  worse or the same, then the node corresponding to the partial solution that created the subproblem is pruned (pruning by bound). A third way of pruning is when a subproblem is solved, i.e., an integral solution with the objective function value  matching 
the optimum LP relaxation value of the subproblem is found
 (pruning by optimality). The LP relaxation of the problem at the root node is solved using the primal simplex algorithm~\cite{Nash}. The LP relaxations of the subproblems 
 created at the non-root nodes 
 are solved using the dual simplex algorithm taking advantage of warm starts~\cite{Banciu2011}.
Every time  B\&B finds a solution with a better objective function value than that of the incumbent best solution, the best solution is updated with the new solution. If the goal is to find an optimum solution, then  B\&B can be  stopped as soon as it finds a solution whose objective function value is equal to the best known lower bound for the optimum value of ILP~(\ref{eqn:geneqILP}) because this solution must be optimum. For more details, see Chapter 7 of~\cite{Wolsey1998}.

 All optimum solutions of ILP~(\ref{eqn:geneqILP}) can be enumerated by using a depth-first search B\&B that branches on the integer values of the variables  if the optimum value of ILP~(\ref{eqn:geneqILP}) is known in advance and say is equal to $z^*$. This is done by pruning a partial solution corresponding to a subproblem if and only if  $z^*$ is strictly smaller than the optimum LP relaxation value of the subproblem or the subproblem is infeasible. When enumerating all optimum solutions no partial solution is pruned by optimality.
 This version of the B\&B algorithm was used in~\cite{Geyer2014} to classify orthogonal arrays  $\oan$ up to isomorphism and a weaker form of isomorphism for many $N,k,s,t$ combinations.
It was also used in~\cite{Bulutoglu2008,Bulutoglu2016,Margot2003b} to classify all non-isomorphic
 $\oan$, covering arrays with the minimum number of rows $\ca$, packing arrays with the maximum number of rows $\pa$,  $4$-$(10,5,1)$-covering designs with the minimum number of sets (blocks), and all $\oant$ up to OD-equivalence for 
 many $N,k,s,t$ combinations. However,  the bottleneck classifications of  all  $\oant$ up to OD-equivalence in~\cite{Bulutoglu2016} required using {\tt nauty}~\cite{McKay2013,McKayP} to remove OD-equivalent $\oant$. In this paper, we also use this version of a B\&B 
algorithm to directly classify all bottleneck $\oant$ in~\cite{Bulutoglu2016} up to OD-equivalence  without resorting to {\tt nauty}~\cite{McKay2013,McKayP} for removing OD-equivalent $\oant$. (The definitions of $\oan$, isomorphism of $\oan$, and OD-equivalence of $\oant$ are deferred until Section~\ref{sec:OA}.) Throughout this paper, when we refer to a B\&B algorithm we mean a depth-first 
search B\&B that branches on the integer values of the variables  targeted to find all optimum solutions of an ILP.

The group of all  permutations of the variables of ILP~(\ref{eqn:geneqILP})
 that map feasible points onto feasible points and preserve the objective function value of each feasible point  is called the {\em  symmetry group} of  
ILP~(\ref{eqn:geneqILP})~\cite{Margot2010,Ostrowski2011}. 
For a subgroup $G$ of the symmetry group of ILP~(\ref{eqn:geneqILP}), two (partial) solutions $\x_1$ and $\x_2$ of  ILP~(\ref{eqn:geneqILP}) are called  {\em isomorphic under the action of $G$} if $g(\x_1)=\x_2$ for some $g \in G$, where $$g((x_1,\ldots,x_n)^{\top})=(
 x_{g^{-1}(1)},\ldots,x_{g^{-1}(n)})^{\top}.$$ Similarly, two subproblems of  ILP~(\ref{eqn:geneqILP}) are called {\em isomorphic subproblems} if they are created from isomorphic partial solutions. Clearly, the LP relaxation of isomorphic subproblems have the same optimum objective function value and feasibility status. Hence,
when the symmetry group of  
ILP~(\ref{eqn:geneqILP}) is large, a B\&B algorithm wastes time by solving 
the LP relaxations of a large number of isomorphic subproblems created from the same number of isomorphic partial solutions. To address this issue,  
 Margot~\cite{Margot2002,Margot2003b,Margot2003a,Margot2007} developed the B\&B with isomorphism pruning algorithm  that finds a set of  all non-isomorphic optimal solutions to  ILP~(\ref{eqn:geneqILP}) by solving the LP relaxation of only the unique subproblem created from the unique  lexicographically minimum partial solution under the action of $G$, where the lexicographical ordering of partial solutions is defined as follows.
  \begin{defn}
 Let $\x$ and $\x'$ be two partial solutions of ILP~(\ref{eqn:geneqILP}) in a B\&B search tree. Let $i_1<\cdots<i_{r_1}$ and $i'_1<\cdots<i'_{r_2}$ be the indices of the variables in $\x$ and $\x'$ that are fixed by branching decisions, and  $\gamma=\min\{r_1,r_2\}$. We say that $\x$ is {\em lexicographically smaller} than $\x'$ if one of the following two conditions is satisfied.
 \begin{enumerate}
\item The first non-zero entry in  $(i_1-i'_1,\ldots,i_{\gamma}-i'_{\gamma})$ is negative.
\item $(i_1,\ldots,i_{\gamma})=(i'_1,\ldots,i'_{\gamma})$ and 
the first non-zero entry in  $(x_{i_1}-x'_{i'_{1}},\ldots,x_{i_{\gamma}}-x'_{i'_{\gamma}})$ is positive.
 \end{enumerate}
 \end{defn} 
  When a B\&B algorithm always selects the minimum index non-fixed variable for branching (called {\em minimum index branching}), then removing a partial solution
 (node) of the B\&B search tree if it is not lexicographically minimum under the action of $G$ results in a B\&B tree whose all feasible leaves are a set of all non-isomorphic optimal solutions~\cite{Margot2007}. 
 This is true because when minimum index branching is implemented, each lexicographically minimum node under the action of $G$ has a unique lexicographically minimum parent node under the same action. Margot~\cite{Margot2007} developed an algorithm based on group theory to decide whether a partial solution is lexicographically minimum under the action of $G$. This algorithm is used within  B\&B with isomorphism pruning  to prune isomorphic partial solutions.
  
 Even when relatively small size groups are used, the search tree of  B\&B with isomorphism pruning is much smaller  than that of  B\&B only, causing huge reductions in computation times~\cite{Margot2007}.  
However, to correctly find or classify solutions that are optimal or prove infeasibility, it is necessary that the symmetry group of
ILP~(\ref{eqn:geneqILP})  or one of its subgroups is used.
Whether  two subproblems are deemed isomorphic depends on the subgroup used within B\&B with isomorphism pruning. Subproblems that are inherently isomorphic may be deemed not to be  isomorphic if a smaller subgroup is used.
Consequently, using a larger subgroup results in a B\&B tree with a smaller number of nodes where LP relaxations must be solved.
Hence, it is desirable to find the symmetry group of a given ILP. If this is not possible,  finding larger subgroups is more desirable (finding the symmetry group of an ILP is an NP-hard problem~\cite{Margot2010}).
One subgroup of the symmetry group of an ILP is the \emph{formulation symmetry group}. Finding this group is as hard as the graph isomorphism
problem, which is not known to be solvable in polynomial time.
 The formulation symmetry group of  ILP~(\ref{eqn:geneqILP}) 
given in~\cite{Arquette2016} is defined to be
\begin{equation}\label{eqn:G}
 G(\A,\bb,\B,\dd,\cc)= \left\{\pi \  | \ \pi(\cc)=\cc, \ \exists \+ \sigma \ \mbox{with}\ \A(\pi, \sigma)=\A,\ \B(\pi, \sigma)=\B,
\   \sigma\left[{\bb\atop \dd}\right]=\left[{\bb\atop \dd}\right]\right\},
\end{equation}
where
$\A(\pi,\sigma)$ is the 
matrix obtained by permuting the columns of $\A$ with $\pi$ followed by a permutation of its rows
 with $\sigma$.
 For a more general definition that covers mixed integer non-linear programs (MINLPs), see~\cite{Liberti2012}. 
 Define the formulation symmetry group of a generic  LP~(\ref{eqn:geneqILP}) to be $G(\A,\bb,\B,\dd,\cc)$, where LP~(\ref{eqn:geneqILP}) is obtained by dropping the $\x \in \mathbb Z^n$ constraint in ILP~(\ref{eqn:geneqILP}).
  Throughout the paper, whenever  ILP~(\ref{eqn:geneqILP}) is referred to as  LP~(\ref{eqn:geneqILP}) it is understood that  LP~(\ref{eqn:geneqILP}) is obtained by dropping the integrality constraints in  ILP~(\ref{eqn:geneqILP}).
If ILP~(\ref{eqn:geneqILP}) or LP~(\ref{eqn:geneqILP}) has no equality constraints, then define its formulation symmetry group 
to be 
\begin{equation}\label{eqn:Gabc}
\Gabc=\left\{\pi \  | \ \pi(\cc)=\cc, \ \exists \+ \sigma \ \mbox{with}
\ \B(\pi, \sigma)=\B,
\   \sigma(\dd)=\dd\right\}.
\end{equation}
 Margot~\cite{Margot2010} and Pfetsch and Rehn~\cite{Rehn2017} described  methods for finding the formulation symmetry group $\Gabc$.
  Each of these methods can be used to find $\GABbcdn$. In~\cite{Liberti2012}, a more general algorithm for finding the formulation symmetry group of an MINLP is described.
 
 Another subgroup of the symmetry group of an ILP is the symmetry group of its LP relaxation, where the two groups may or may not be the same. The subgroup property follows directly from the  following definition and the definition of the symmetry group of an ILP.
\begin{defn}
 Let $\mathcal{F}$ be the feasible set of  LP~(\ref{eqn:geneqILP}) and
\begin{align*}
\Glp=&\ \{\pi \in S_n \  | \ \pi(\x) \in \mathcal{F} \   \text{and }\cc^{\top}\pi(\x)=\cc^{\top}\x\ \text{ for all } \x \in \mathcal{F}\},\\
\Glp_{\mathcal{F}}=&\ \{\pi \in S_n \  | \ \pi(\x) \in \mathcal{F} \  \text{ for all } \x \in \mathcal{F}\},
\end{align*} 
where $S_n$ is the set of all permutations of indices $\{1,\ldots,n\}$.
Then the group $\Glp$ is called the {\it symmetry group} of  LP~(\ref{eqn:geneqILP}), and $\Glp_{\mathcal{F}}$ is called the {\it symmetry group of the feasible set} of LP~(\ref{eqn:geneqILP}).
\end{defn}
 Hence, 
$\Glp$ of an LP is completely determined by its feasible set and its objective function. In particular, 
$\Glp$ of an infeasible LP with $n$ variables is $S_n$.  Clearly, the formulation symmetry group $\GABbcdn$ is a subgroup of $\Glp$ of LP~(\ref{eqn:geneqILP}), and
 $\Glp$ of LP~(\ref{eqn:geneqILP}) is a subgroup of the symmetry group of ILP~(\ref{eqn:geneqILP}).
This makes it viable to use $\Glp$ or $\GABbcdn$ within B\&B with isomorphism pruning to find a set of all non-isomorphic solutions to ILP~(\ref{eqn:geneqILP}). However, $\GABbcdn$ may be too small to reap the full benefits of B\&B with isomorphism pruning. Hence, it is essential to develop methods that find $\Glp$ in general.

In Section~\ref{sec:orbitproj}, we prove that $\Glp$ of an LP
coincides with its formulation symmetry group  if the feasible set of the LP is full dimensional and it has no redundant inequalities. Therefore, the method in~\cite{Margot2010,Rehn2017}  can be
 used to find the symmetry group of a full-dimensional LP after removing all redundant inequalities. Conversely, the formulation symmetry group of an LP can miss inherent symmetries if the LP has redundant constraints. This is discussed in Section~\ref{sec:express}.

Different LP formulations based on the same variables can have the same feasible set. When this happens we say that the two LP formulations \emph {define} the same feasible set.
We say an LP in the form of LP~(\ref{eqn:geneqILP}) is in  \emph{standard form} if it is feasible, has no redundant  constraints, 
and none of the inequalities in $\B\x\leq \dd$ is  satisfied by every feasible $\x$ as an equality. 
Section~\ref{sec:express} describes a method for defining the feasible set of a given feasible LP by an LP in  standard form and with the same objective function.
There is no known general method for finding the symmetry group of a feasible LP that is not full dimensional. 
In Section~\ref{sec:orbitproj}, we describe a method based on orthogonal projection matrices that finds the symmetry group of a non-full-dimensional LP in standard form.

 In Section~\ref{sec:OA}, we 
define orthogonal arrays (OAs) and describe the isomorphism and OD-equivalence operations that map OAs to OAs. 
In Section~\ref{sec:OAapply}, we analytically characterize a subgroup of the
LP relaxation symmetry group $\Glp$ of an OA defining ILP   in terms of the isomorphism and OD-equivalence operations.
In Section~\ref{sec:comparisons}, we apply the Section~\ref{sec:orbitproj} method to compute the  LP relaxation symmetry groups of many cases of
 an OA defining ILP formulation from~\cite{Bulutoglu2008}. There is an OA defining ILP formulation in~\cite{Bulutoglu2016} with the objective function $0$ and without redundant constraints. We then make speed comparisons between using $\Glp$ with this ILP formulation from~\cite{Bulutoglu2016} 
  and two other group/formulation combinations from~\cite{Geyer2014} within  the B\&B with isomorphism pruning algorithm from~\cite{Margot2007} for enumerating OAs up to OD-equivalence, isomorphism, and a weaker form of isomorphism. 
In particular, our method reduced the computation time  to find all OD-equivalence 
classes of OA$(160,8,2,4)$ and OA$(176,8,2,4)$ by  factors of $1/(2.16)$ and $1/(1.36)$ compared to the fastest known method in~\cite{Bulutoglu2016}. These are the largest $2$-symbol, strength $4$ cases for which classification results are available and yet only symmetry exploiting methods have successfully generated them.
Moreover, for most OA defining ILPs with only inequalities 
that we considered,  speedups
 gleaned from exploiting the additional LP relaxation symmetry captured by adding slack variables drastically overcome the additional computational burdens due to the added variables.
 In Section~\ref{sec:future}, we discuss the major findings of this paper and propose a future research project. 
 
 Throughout the paper, $S_n$ is either the symmetric group of degree $n$ or an isomorphic copy of it. If the action of $S_n$ is not defined within a paragraph, then it can be assumed that $S_n$ is the abstract symmetric group of degree $n$ within that paragraph.
\section{A method for putting a feasible LP in standard form}\label{sec:express}
In this section, we provide  a method for putting a  feasible  LP in standard form. First, we need the well-known  Theorem \ref{thm:express}. Theorem~\ref{thm:express} leads to Method \ref{meth:equality} for finding all equality constraints of a
feasible LP. 

\begin{thm}\label{thm:express}
Let $P\neq \emptyset$, $P \subseteq \mathbb{R}^n$ be the feasible set of a system
of constraints
\begin{equation}\label{eqn:cons}
\begin{array}{l}
   \A \x=\bb, \\
   \B \x \leq \dd
\end{array}
\end{equation}
and $\boldsymbol{\beta}_i^{\top}$ be the $i$th row of $\B$ for $i=1,\dots,m$.
Then $P$ is full dimensional in the affine space $\A \x=\bb$ if and only if there is a sequence of feasible points $\{\x_i\}_{i=1}^m$ of constraints~(\ref{eqn:cons}) such that
$\boldsymbol{\beta}_i^{\top}\x_i<d_i$.
\end{thm}
\begin{algorithm}[t]
\centering
\caption{\label{meth:equality} Finding all equality constraints of a feasible LP of form~(\ref{eqn:geneqILP})}
\begin{algorithmic}[1]
\State {\bf Input} a feasible LP $L$ of form (\ref{eqn:geneqILP}) with $m\times n$ inequality constraint matrix $\B$.
\For{$i:= 1$ {\bf to} $m$ {\bf step} $1$}
\State {\bf Set} $\boldsymbol{\beta}^{\top}:=\B_i;$ \Comment{$\B_i$ is the $i$th row of $\B$.}
\State {\bf Solve} LP 
\begin{equation*}\label{eqn:lp2}
\begin{array}{ccl}
y_i&:=&\min\limits_{\x}   \boldsymbol{\beta}^{\top}\x  \\
&&\mbox{s.t.}\   \A \x=\bb, \ \B\x\leq \ddd;
\end{array}
\end{equation*}
\EndFor
\For{$i:= m$ {\bf to} 1 {\bf step} $-1$} 
\If{$y_i=d_i$}\Comment{Change the $i$th inequality constraint of $L$ to an equality constraint.} \label{step:yidi}
\State  {\bf Append} 
\begin{equation*}
\A := \left[ \begin{array}{c}\A\\  \boldsymbol{\beta}^{\top} \end{array}\right]\text{~~and~~}\bb:=\left(\begin{array}{c}\bb\\  d_i \end{array} \right);
\end{equation*}
\State {\bf Delete} the $i$th row of $\B$ and the $i$th entry of $\ddd$;
\EndIf 
\EndFor
\State {\bf Output} $L$.
\end{algorithmic}
\end{algorithm}
\begin{rem}
For practical purposes,  $y_i\in[d_i-10^{-6},d_i+10^{-6}]$ can be used instead of $y_i=d_i$ in Step~\ref{step:yidi} of Method~\ref{meth:equality}.
\end{rem}
 It is always possible to inscribe a highly-symmetric polytope inside an asymmetric polytope so that the formulation symmetry group of the resulting system of constraints is much smaller. This idea is formalized in the following theorem. We skip the proof of this well-known result.
\begin{thm}\label{thm:redundant}
The formulation symmetry group of every bounded LP $L$ with a finite number of constraints can be reduced to the identity permutation  by adding redundant inequalities.
\end{thm}
By Theorem~\ref{thm:redundant}, redundant constraints can mask inherent symmetries of an LP. Hence, it is essential to  remove the redundant inequalities before computing the formulation symmetry group. LPs in standard form have no redundant constraints.  
Method \ref{meth:express} finds an LP in standard form that defines the feasible set of a given feasible LP having the same objective function.

\begin{algorithm}[t]
\centering
\caption{\label{meth:express} Putting a feasible LP $L$ of form~(\ref{eqn:geneqILP}) in standard form}
\begin{algorithmic}[1]
\State {\bf Input} a feasible LP $L$ of form (\ref{eqn:geneqILP}).
\State {\bf Apply} Method~\ref{meth:equality} to $L$ and overwrite $L$ with the result;
\State {\bf Remove} all redundant inequality constraints from $L$ by solving a sequence of LPs;
\State {\bf Remove} a set of all redundant equality constraints from $L$ by using Gaussian elimination;
\State {\bf Output} $L$.
\end{algorithmic}
\end{algorithm}
\section{A method for finding the symmetry group of a feasible LP}
\label{sec:orbitproj}
The symmetry group $\Glp$ of an LP is completely determined by its feasible set and objective function. Then the symmetry group $\Glp$ of a given feasible LP can be found by finding the symmetry group  of an LP in standard form that has the same feasible set and objective function as the given LP. Such an LP in standard form can be obtained by applying Method~\ref{meth:express} from Section~\ref{sec:express}.
 Next, we describe a method for finding the symmetry group $\Glp$ of LP~(\ref{eqn:geneqILP}) in standard form. 
  Let Row$(\A)$ be
   the row space of $\A$ and
\begin{equation}\label{eqn:projmatrix1}
\PP_{\A^{\top}}=\A^{\top}(\A\A^{\top})^{-1}\A
\end{equation}  
  be the orthogonal projection matrix onto Row$(\A)$. 
Let  $p$ be the number of rows of $\A$. Thus,  
$p={\rm rank}(\A)$, i.e., $\A$ has $p$ linearly independent rows. For a vector $\boldsymbol{\vv} \in \mathbb{R}^p$, let $\mbox{diag}(\boldsymbol{\vv})$ be the diagonal matrix whose $i$th diagonal entry is $v_i$ for $i \in \{1,\ldots,p\}$. Let $\boldsymbol{\sigma}$  be a vector of singular values of $\A$ such that 
$\A=\U\D\V^{\top}$ is a singular value decomposition of $\A$, where 
$\U\U^{\top}=\U^{\top}\U=\I_p$,  $\V\V^{\top}=\V^{\top}\V=\I_n$, and the $p \times n$ matrix
$$\D=\begin{bmatrix}[l]\mbox{diag}(\boldsymbol{\sigma}) & \0 
 \end{bmatrix}$$ 
 is based on the all zeros matrix  $\0$ of appropriate dimension~\cite{Leon2002}. Then, 
equation~(\ref{eqn:projmatrix1}) simplifies to
\begin{align}\label{eqn:projmatrix2}
 \PP_{\A^{\top}}&=\V\D^{\top}\U^{\top}(\U\D \V^{\top}\V\D^{\top}\U^{\top})^{-1}
 \U\D\V^{\top} \nonumber \\
 &=\V\D^{\top}(\D\D^{\top})^{-1}\D\V^{\top}=\V\I_n^{(p)}\V^{\top},
\end{align}
where 
\begin{equation*}
\I_n^{(p)} = \begin{bmatrix}[l]\I_{p \times p} & \0 \\ \0 & \0 \end{bmatrix}
\end{equation*}
is $n \times n$. 
Equation~(\ref{eqn:projmatrix2}) should be used to compute $\PP_{\A^{\top}}$ as it does not involve matrix inversion, 
 leading to improved accuracy especially when $\A$ 
 is ill-conditioned.

Let $S_n$ be the group of all permutations of coordinates of column vectors in $\mathbb{R}^n$.
Observe that each $\pi \in S_n$ is a linear transformation from $\mathbb{R}^n$ to $\mathbb{R}^n$. Let $\ppi$ be the matrix of $\pi\in S_n$ with respect to the standard basis $\{\e_1,\ldots,\e_n\}$. 
Since $(\ppi\vv)^{\top}=\vv^{\top}\ppi^{\top}$, right multiplication of $\vv^{\top}$ by $\ppi^{\top}$  permutes the coordinates 
of the row  vector $\vv^{\top}$.
 The {\it automorphism group} of an $n \times n$ matrix $\M$, denoted by $G_ {\M}$, is the set of all $\pi\in S_{n}$ that send $\M$ to itself when the rows and the columns of $\M$ are permuted according to $\pi$. So,
$$G_ {\M}= \{\pi\in S_n\  |\  \ppi\M\ppi^{\top}=\M\}.$$
For a vector space $V \subseteq \mathbb{R}^n$, define Stab$(V)=\{\pi \in S_n\  |\  
\ppi \vv \in V\ \forall\  \vv \in V\}$.
Then we have the following lemma. 
\begin{lem}\label{lem:stabrowA}
Let $\A$ be an $m\times n$ matrix with full row rank and  $\PP_{\A^{\top}}$ be the orthogonal projection 
matrix onto Row$(\A)$. 
 Then $G_ {\PP_{\A^{\top}}}={\rm Stab}({\rm Row}(\A))$.
\end{lem}
\begin{pf}
To prove ${\rm Stab}({\rm Row}(\A)) \subseteq 
G_ {\PP_{\A^{\top}}}$, let $\pi \in {\rm Stab}({\rm Row}(\A))$.
Then, since $\pi \in {\rm Stab}({\rm Row}(\A))$ and $\ppi$ is an  invertible matrix, ${\rm Row}(\A)={\rm Row}(\A\ppi^{\top})$. Hence, the set of  rows of 
$\A\ppi^{\top}$ is a basis for ${\rm Row}(\A)$. Moreover,  $\ppi^{\top}\ppi=\ppi\ppi^{\top}=\I$ as 
every permutation matrix  is an orthogonal matrix.
Then, $$\PP_{\A^{\top}}=\PP_{\ppi\A^{\top}}=(\A\ppi^{\top})^{\top}(\A\ppi^{\top}(\A\ppi^{\top})^{\top})^{-1}\A\ppi^{\top}=\ppi\A^{\top}(\A\ppi^{\top}\ppi\A^{\top})^{-1}\A\ppi^{\top}.$$ 
Hence, $$\PP_{\A^{\top}}=\PP_{\ppi\A^{\top}}=\ppi\A^{\top}(\A\A^{\top})^{-1}\A\ppi^{\top}=\ppi\PP_{\A^{\top}}\ppi^{\top},$$
and $\pi \in G_ {\PP_{\A^{\top}}}$.

To prove $G_ {\PP_{\A^{\top}}}\subseteq {\rm Stab}({\rm Row}(\A))$,
let $\pi \in G_ {\PP_{\A^{\top}}}$. Then 
\begin{align}\label{eqn:PP}
\PP_{\A^{\top}}=&\ \ppi\PP_{\A^{\top}}\ppi^{\top} \nonumber\\
\PP_{\A^{\top}}\ppi=&\ \ppi\PP_{\A^{\top}}.
\end{align}
Let  $ {\rm Col}(\M)$
 of a matrix $\M$ be the column space of $\M$ and $\w\in {\rm Row}(\A)$, where $\w$ is written as a column vector. Then $\w=\PP_{\A^{\top}}\w$, and by~(\ref{eqn:PP}), we have
$$
\ppi\w=\ppi \PP_{\A^{\top}}\w=\PP_{\A^{\top}}\ppi\w.
$$
Hence, $\ppi\w \in {\rm Col}(\PP_{\A^{\top}})={\rm Col}(\A^{\top})={\rm Row}(\A)$, and
 $\pi \in {\rm Stab}({\rm Row}(\A))$. 
\qed
\end{pf}

Let $G(\B,\ddd,\cc)$ be the formulation symmetry group as defined in~(\ref{eqn:Gabc}) and $$G(\B,\ddd,\cc)=\{\pi\in S_n\  |\  \pi(\cc)=\cc\}$$ when $\B$ is the empty matrix and $\ddd$ is the empty vector.
 Let $\GABocd$  be the largest subgroup of  $G_{\PP_{\A^{\top}}}$  that preserves $\cc$ and the set of  inequalities in  $\B\x\leq\dd$. 
Then by Lemma~\ref{lem:stabrowA},
\begin{equation}\label{eqn:intersection}
\GABocd=G_ {\PP_{\A^{\top}}}\bigcap G(\B,\ddd,\cc)=
{\rm Stab}({\rm Row}(\A))\bigcap G(\B,\ddd,\cc).
\end{equation}

An {\em automorphism}  of a vertex colored, edge colored graph is a permutation of its vertices that maps 
adjacent vertices to adjacent vertices and preserves  vertex and edge colors. The set of all such permutations forms a group called the {\em automorphism group} of the graph. 
 Method \ref{meth:GPAT} computes $\GABocd$ as the intersection of the automorphism group of a vertex colored, edge colored  graph and $G(\B,\ddd,\cc)$. The formulation symmetry group $G(\B,\ddd,\cc)$ can  be computed as the automorphism group of a vertex colored, edge colored 
 graph with $n+m'$ vertices, where $\B$ is $m' \times n$~\cite{Margot2010, Rehn2017}.
 Edge coloring in this computation is necessary only if the number of distinct numerical values in the matrix $\B$ is larger than two~\cite{Margot2010}.
 
\begin{algorithm}[t]
\centering
\caption{\label{meth:GPAT} Computing $\GABocd$ of an LP $L$ of form~(\ref{eqn:geneqILP}) in standard form}
\begin{algorithmic}[1]
\State {\bf Input} $\A,\B,\ddd,\cc$ from an  LP $L$ of form~(\ref{eqn:geneqILP}) in standard form.
\State {\bf Compute} a singular value decomposition $\A=\U\D\V^{\top}$ and $\PP_{\A^{\top}}=\V\I_n^{(p)}\V^{\top}$; \label{step:svd}

\State {\bf Label} each distinct numerical value in $\PP_{\A^{\top}}$ with a distinct color;\label{step:labeledge} 
\State {\bf Set} $nce:=$ number of distinct colors in Step~\ref{step:labeledge};\label{step:nce}
\State {\bf Initialize} $\mathcal{G}(\PP_{\A^{\top}})$  to be the graph with $n$ vertices and no edges;
\For{$i:= 1$ {\bf to} $(n-1)$ {\bf step} $1$}\label{step:edgestart}
\For{$j:= (i+1)$ {\bf to} $n$ {\bf step} $1$}
\For{$\ell:= 1$ {\bf to} $nce$ {\bf step} $1$}
\If{the $(i,j)$th entry of $\PP_{\A^{\top}}$ is labeled with color $\ell$}
\State {\bf Put} an edge between $i$th and $j$th vertices of $\mathcal{G}(\PP_{\A^{\top}})$ with color $\ell$; \label{step:put}
\EndIf
\EndFor
\EndFor
\EndFor \label{step:edgend}
\State {\bf Label} each distinct numerical value in $\cc$ with a distinct color;\label{step:labelvertex}
\State {\bf Set} $ncv:=$ number of distinct colors in Step~\ref{step:labelvertex};\label{step:ncv}
\For{$i:= 1$ {\bf to} $n$ {\bf step} $1$}\label{step:vertexstart}
\For{$\ell:= 1$ {\bf to} $ncv$ {\bf step} $1$}
\If{the $i$th entry of $\cc$ is labeled with color $\ell$}
\State {\bf Color} vertex $i$ of $\mathcal{G}(\PP_{\A^{\top}})$ with color $\ell$;
\EndIf
\EndFor
\EndFor\label{step:vertexend}
\State {\bf Compute} the automorphism group $H_{\PP_{\A^{\top}}}$ of the vertex colored, edge colored graph $\mathcal{G}(\PP_{\A^{\top}})$;\label{step:edgecolored} 
\State {\bf Compute} $G(\B,\ddd,\cc)$ by computing  the automorphism group of a graph~\cite{Margot2010,Rehn2017};\label{step:formulation}
\State \label{step:intersect}{\bf Compute} $\GABocd:=
H_{\PP_{\A^{\top}}}\bigcap G(\B,\ddd,\cc)$;
\State {\bf Output} $\GABocd$.
\end{algorithmic}
\end{algorithm}

Given two graphs, the {\em graph isomorphism} (GI) problem asks whether one can be obtained from the other by permuting vertices.
Finding the generators of the automorphism group of a graph is known to be equivalent to the GI problem~\cite{Luks}. Finding the generators of the
intersection of two groups is also  equivalent to the GI problem~\cite{Comp}. 
 It is easy to see that the GI problem is in NP. On the other hand,  it is not known whether the GI problem is NP-complete. It is also not known if the GI problem is in P.
All the known algorithms for the GI problem have 
 exponential worst-case running times.
 For algebraic techniques that compute the generators for the automorphism group of a graph, see~\cite{McKay1981,McKayP}, and for the intersection of two groups,  see~\cite{Rehn2010}.

There is available software that can be used in implementing Method~\ref{meth:GPAT}.
In Step~\ref{step:put},  edge coloring can be implemented
by using a vertex colored graph with $n\lceil\log_2{(nce+1)}\rceil$ vertices, where $nce$ is the number of distinct numerical values in $\PP_{\A^{\top}}$ and $n$ is the number of columns  of the square matrix $\PP_{\A^{\top}}$~\cite{McKay2013}.  The subgroup $H_{\PP_{\A^{\top}}}$ of the automorphism group $G_ {\PP_{\A^{\top}}}$ of  $\PP_{\A^{\top}}$ that preserves $\cc$ in Step~\ref{step:edgecolored}
and the formulation symmetry group $G(\B,\ddd,\cc)$ in Step~\ref{step:formulation} can both be computed as the automorphism groups of their corresponding vertex colored, edge colored graphs by using \verb!nauty!~\cite{McKay2013, McKayP}. 
In Step~\ref{step:intersect}, 
the intersection  can be computed  by using \verb!GAP!~\cite{GAP2013}.

\begin{defn}\label{def:GAbBdc}
For an LP $L$ of form~(\ref{eqn:geneqILP}) in standard form, $\Glpp$ is defined to be the largest subgroup of $\Glp$ of $L$ that preserves the vector $\cc$.
\end{defn}
\begin{lem}\label{lem:subset}
For an LP $L$ of form~(\ref{eqn:geneqILP}) in standard form, 
$\Glpp\leq\GABocd$.
\end{lem}
\begin{pf}
By the definition of standard form, $L$ is feasible. 
Let $\pi \in \Glpp$ and $\A$ be an $m \times n$ matrix. Then each such $\pi$
must preserve the facets $\B\x\leq\ddd$ and the objective function coefficient vector $\cc$. An orthogonal linear transformation preserves ${\rm Row}(\A)$ if and only if it preserves ${\rm Null}(\A)={\rm Row}(\A)^{\perp}$.
Since $\pi$ is an orthogonal linear transformation, by equations~(\ref{eqn:intersection}), it suffices to prove that  $\pi$ preserves  ${\rm Null}(\A)$. 
Since the feasible set of $L$ is a full-dimensional 
polytope in an affine space of dimension $n-m$, the feasible set of 
$L$ contains $n-m+1$ affinely independent points $\x_j\in \mathbb{R}^n$ for $j=0,\ldots,n-m.$ Then the vectors $\vv_j=\x_j-\x_0 \in {\rm Null}(\A)$
for $j=1,\ldots,n-m$  are linearly independent. Moreover, 
$${\rm Span}(\vv_1,\ldots,\vv_{n-m})={\rm Null}(\A)$$ as  ${\rm dim(Null}(\A))=n-m$ by the rank-nullity theorem~\cite{Leon2002}, where  ${\rm dim}(V)$ denotes the dimension
of a vector space $V$.

 Since $\pi$ is a linear transformation from $\mathbb{R}^n$ to $\mathbb{R}^n$, 
\begin{equation}\label{eqn:null}
\pi({\rm Span}(\vv_1,\ldots,\vv_{n-m}))={\rm Span}(\pi(\vv_1),\ldots,\pi(\vv_{n-m}))=\pi({\rm Null}(\A)),
\end{equation}
where for a vector space $V\subseteq \mathbb{R}^n$ and a linear transformation $T$ from $\mathbb{R}^n$ to $\mathbb{R}^n$
$$T(V)=\{\w\in \mathbb{R}^n \ |\  \w=T(\vv) \text{ for some } \vv \in V\}. $$
Observe that $\pi(\vv_j)=\pi(\x_j)-\pi(\x_0)$ and $\A \pi(\vv_j)=\A\pi(\x_j)-\A\pi(\x_0)=\bb-\bb=\0.$ Then, 
$\pi(\vv_j) \in {\rm Null}(\A)$ for $j =1,\ldots,n-m$ and by equations~(\ref{eqn:null})
 $$\pi({\rm Null}(\A)) \subseteq {\rm Null}(\A).$$
Since $\pi$ is an invertible linear transformation and $\vv_1,\ldots,\vv_{n-m}$ are linearly independent, the vectors $\pi(\vv_1), \ldots, \pi(\vv_{n-m})$ are linearly independent. Consequently, $${\rm dim}(\pi({\rm Null}(\A)))={\rm dim}({\rm Span}(\pi(\vv_1),\ldots,\pi(\vv_{n-m})))=n-m.$$
Hence, since ${\rm dim(Null}(\A))=n-m$, 
$$\pi({\rm Null}(\A))= {\rm Null}(\A).$$ \qed
\end{pf}
Let $H\leq \GABocd$,
 $\mathcal{B}=\{\e_1,\dots, \e_{n}\}$ be the standard  basis of $\mathbb{R}^{n}$, and $O_1,\dots, O_r$ be the orbits in $\mathcal{B}$ under the action of $H$, i.e., for each $\x, \y \in \mathcal{B}$ there exists  $g \in H$ such that $\x=g(\y)$ if and only if $\x,\y \in O_i$ for some $i$.  The fixed subspace of $\mathbb{R}^{n}$ under the action of $H$ is defined as
\begin{equation*}\label{eqn:Fixnew}
{\rm Fix}_{H}(\mathbb{R}^{n}):=\{\x \in \mathbb{R}^{n}\ |\ \gamma\x=\x\text{ for all $\gamma
 \in H$}\}.
\end{equation*}
Lemma 3 in~\cite{Bodi2013} implies that
\begin{equation}\label{eqn:Fix}
{\rm Fix}_{H}(\mathbb{R}^{n})={\rm Span}(\beta(O_1),\dots,\beta(O_r)),
\end{equation}
where for a set $S$ of vectors
\begin{equation}\label{eqn:beta}
\beta(S)=\frac{\sum_{\vv\in S}\vv}{|S|}.
\end{equation}
Let $\E$ be the orthogonal projection matrix onto ${\rm Span} (\beta(O_1),\dots,\beta(O_r))$ with respect to $\mathcal{B}$.
Then
\begin{eqnarray}\label{eqn:span}
E_{ij}=\left\{
\begin{array}{cl}
 \frac{1}{|O_{i,j}|}& \quad \text{if $i$ and $j$ belong to the same orbit $O_{i,j}\in \{O_1,\dots, O_r\}$,}\\
  0&\quad  \text{otherwise.}
\end{array}\right.
\end{eqnarray}
The matrix $\E$ uniquely identifies ${\rm Fix}_{H}(\mathbb{R}^{n})$. Let $L$ be an LP of form~(\ref{eqn:geneqILP}) in standard form, $\mathcal{F}^{L}$ be its feasible set, and $\mathcal{T}^{L}_{{\rm Fix}_{H}}=\mathcal{F}^{L}\cap {\rm Fix}_{H}(\mathbb{R}^{n}).$ Then
\begin{equation*}\label{eqn:TFix}
\mathcal{T}^{L}_{{\rm Fix}_{H}}=\{\x \in {\mathbb R}^n\  | \ \left(\I- \E\right) \x=\0 \text{ and $\x$ is a feasible point of } L\}.
 \end{equation*}
 Now, we have the following theorem.
\begin{thm}\label{thm:thm}
 Let $L$ be an LP of form~(\ref{eqn:geneqILP}) in standard form and $H\leq \GABocd$. Then 
 $ \mathcal{T}^{L}_{{\rm Fix}_{H}}$ is non-empty if and only if 
 $H\leq\Glpp$.
\end{thm}
\begin{pf}
By the definition of standard form, $L$ is feasible. Let $\vv_0 \in  \mathcal{T}^{L}_{{\rm Fix}_{H}}$, let $$\mathcal{F}_{\B\x\leq \dd}=\{\x\in \mathbb{R}^n\ |\ \B\x\leq 
\ddd\},$$ and for a set $\mathcal{S} \subseteq \mathbb{R}^n$ and a vector $\uuu\in \mathbb{R}^n$, let $ \mathcal{S}+\uuu=\{\vv\in \mathbb{R}^n\ |\  \vv=\s+\uuu\ \text{for some}\ \s\in \mathcal{S}\}.$ Since $\vv_0$   is  in the feasible set of $L$, $\x \in \mathbb{R}^n$ is in the feasible set of $L$ if and only if
\begin{equation*}\label{eqn:particular}
\x=\vv_0+\vv
\end{equation*}
for some $\vv \in {\rm Null}(\A) \bigcap \left(\mathcal{F}_{\B\x\leq\ddd}-\vv_0\right)$.
  Let $\x$ be in the feasible set of $L$ and $\pi\in H.$
  Each  $\pi \in H\leq\GABocd$ preserves $\cc$ and  $ \mathcal{F}_{\B\x\leq\ddd}-\vv_0$ as it preserves $\cc$, $ \mathcal{F}_{\B\x\leq\ddd}$, and $ \vv_0$. By equations~(\ref{eqn:intersection}), $\pi$ preserves  ${\rm Row}(\A)$. Then, ${\rm Null}(\A)={\rm Row}(\A)^{\perp}$ implies that $\pi$ also preserves  ${\rm Null}(\A)$. Hence, $\pi$ preserves $\cc$ and ${\rm Null}(\A)\bigcap \left(\mathcal{F}_{\B\x\leq\ddd}-\vv_0\right).$ Then
  $$ \A\pi(\x)=\A\pi(\vv_0)+\A\pi(\vv)=\A\vv_0=\bbb  $$
  and
  $$\pi(\x)=\pi(\vv_0)+\pi(\vv)=\vv_0+\pi(\vv) \in \mathcal{F}_{\B\x\leq\ddd}.$$
Hence, $\pi \in \Glpp.$ This proves  $H \leq \Glpp$. 

To prove the converse, let $H\leq\Glpp$ and $\x_0$ be a feasible point in  $L$. Let $$O_{\x_0}=\{\y\in\mathbb{R}^n\ |\ \y =h(\x_0) \text{ for some } h\in H \}$$ be the orbit of 
$\x_0$ under the action of $H$ on $\mathbb{R}^n$ and $\beta$ be the orthogonal projection operator onto ${\rm Fix}_{H}(\mathbb{R}^{n})$ as defined  in equations~(\ref{eqn:Fix}) and~(\ref{eqn:beta}). 
Now, since   $\beta(O_{\x_0})$ is a convex combination of feasible points of $L$, $\beta(O_{\x_0})$ is feasible. Hence,  $\beta(O_{\x_0})\in \mathcal{T}^{L}_{{\rm Fix}_{H}}$.\qed
\end{pf}
\begin{cor}\label{cor:T}
  Let  $L$ be an LP of form~(\ref{eqn:geneqILP}) in standard form. Then $ \mathcal{T}^{L}_{{\rm Fix}_{\GABocd}}$ is non-empty if and only if 
 $\GABocd=\Glpp$.
\end{cor}
\begin{pf}
The result follows from Lemma~\ref{lem:subset} and Theorem~\ref{thm:thm}.\qed
\end{pf}
 Method \ref{meth:Glpp} uses the formulation symmetry group  
 $\GABbcdn$ defined  in equation~(\ref{eqn:G}) and the output 
 $\GABocd$ from Method~\ref{meth:GPAT} 
 to find the $\Glpp$ of an LP $L$ of form~(\ref{eqn:geneqILP}) in standard form. Let $\A$ be $m \times n$ and $\B$ be $m' \times n$. Then the formulation symmetry group $\GABbcdn$ can be computed as the automorphism group of a vertex colored, edge colored 
 graph with $n+m+m'$ vertices, where edge coloring is necessary only if the number of distinct numerical values in the matrix $[\A\, \B]$ is larger than two~\cite{Margot2010,Rehn2017}.
  Method~\ref{meth:Glpp} requires 
  finding a double coset 
decomposition of $\GABocd$  with respect to its subgroup 
$\GABbcdn$ and solving either $1$  or $q$ LPs, where $q$ is the 
number of double cosets.
 In terms of computational complexity,  it is not known whether there is a polynomial time algorithm for determining the number of double cosets in a double coset decomposition of a permutation group~\cite{Hoffman}.
 Moreover, the double coset membership problem (i.e., the problem of determining whether a given permutation in a permutation group is in a given double coset) is at least as difficult
 as the GI problem~\cite{Hoffman}. 
 All the known algorithms for computing a double coset decomposition have exponential worst-case running times. 
Method~\ref{meth:Glpp} also requires 
computing the orbits in $\mathcal{B}=\{\e_1,\ldots,\e_n\}$ under the action of $\GABocd$  
in Step~\ref{step:calculate0} and  $G_{\rm ext}$ in Step~\ref{step:calculate}. Given a set $S$ of generators  for a group $G$ acting on a set $\Omega$,
 the orbit $\omega^G$ of an element $\omega \in \Omega$ can be computed in $O(|S||\omega^G|)$ time~\cite{Rehn2010}, where $$\omega^G=\{\omega' \in \Omega\ |\ \omega'=g\omega \text{ for some } g\in G\}.$$  
Using this result, it is easy to see that   the orbits in $\mathcal{B}$ under the action of $\GABocd$ can be computed in $O(|S|n+n^2)$ time.
All computations in Method~\ref{meth:Glpp} involving a group can be implemented in \texttt{GAP}~\cite{GAP2013}, and the feasibility of LPs can be determined by using the primal or dual simplex algorithm implementation in \texttt{CPLEX}~\cite{Cplex}. The following theorem validates Method~\ref{meth:Glpp}.
\begin{algorithm}[t]
\centering
\caption{\label{meth:Glpp} Computing $\Glpp$ of an LP $L$ of form~(\ref{eqn:geneqILP}) in standard form}
\begin{algorithmic}[1]
\State {\bf Input} $\A,\B,\bb,\ddd,\cc$ from an  LP $L$ of form~(\ref{eqn:geneqILP}) in standard form.
\State {\bf Initialize} $i:= 1$;
\State {\bf Compute} $\GABbcdn$ by computing  the automorphism group of a graph~\cite{Margot2010,Rehn2017}; 
\State {\bf Compute} $\GABocd$ by Method~\ref{meth:GPAT};
\State {\bf Compute} the orbits in $\mathcal{B}=\{\e_1,\ldots,\e_n\}$ under $\GABocd$ ;\label{step:calculate0}\Comment{$\mathcal{B}$ is the standard basis.}
\State {\bf Set} $\E$ to be as in equation~(\ref{eqn:span}) for $H=\GABocd$;
\State {\bf Solve}\label{step:determine}  the LP  obtained by adding the constraint $(\I-\E)\x=\0$ to $L$;
\If{the LP in Step~\ref{step:determine} is feasible} 
\State {\bf Set} $\Glpp:=\GABocd$ and {\bf GOTO} Step \ref{meth:GlppOut}; 
\Else \Comment{Compute a $\GABbcdn$-double coset decomposition of $\GABocd$.}
\State \label{step:coset} {\bf Compute} $q$ and a set $\{g_1,\ldots,g_q\}$ so that  
$$
\GABocd=\bigcup_{j=1}^q \GABbcdn g_j \GABbcdn;
$$
\EndIf
\State {\bf Set}  $G_{\rm ext}:= \langle g_1, \ \GABbcdn \rangle$; \Comment{The group generated by $g_1$ and $\GABbcdn$.} \label{step:Gext1}
\State {\bf Compute} the orbits in $\mathcal{B}=\{\e_1,\ldots,\e_n\}$ under $G_{\rm ext}$; \label{step:calculate} 
\State {\bf Set} $\E$  to be as in equation~(\ref{eqn:span}) for $H=G_{\rm ext}$; 
\State {\bf Solve} the LP obtained by adding the constraint $(\I-\E)\x=\0$ to $L$;\label{step:feasible}
\If{the LP in Step~\ref{step:feasible} is feasible} 
\State {\bf Update} $\GABbcdn:=G_{\rm ext}$; 
\EndIf
\State \label{step:return} {\bf Increment} $i:=i+1$;
\If{$i=q+1$} 
\State {\bf Set} $\Glpp:=\GABbcdn$ and {\bf GOTO} Step \ref{meth:GlppOut}; 
\Else \State {\bf Set} $G_{\rm ext}:=\langle g_i,~G(\A,\bb,\B,$ $\ddd,\cc) \rangle$ and {\bf GOTO} Step~\ref{step:calculate}; \label{step:Gext2}
\EndIf
\State {\bf Output} $\Glpp$.\label{meth:GlppOut}
\end{algorithmic}
\end{algorithm}
\begin{thm}
The output of Method~\ref{meth:Glpp} is $\Glpp$.
\end{thm}
\begin{pf}
The set $\mathcal{T}^L_{{\rm Fix}_{\GABocd}}\neq\emptyset$ if and only  if the LP in Step~\ref{step:determine} is feasible. If $\mathcal{T}^L_{{\rm Fix}_{\GABocd}}\neq\emptyset$, then
$\Glpp=\GABocd$ by Corollary~\ref{cor:T}.
  If $\mathcal{T}^{L}_{{\rm Fix}_{\GABocd}}=\emptyset$,  then $$\GABbcdn\leq\Glpp < \GABocd$$ by Lemma~\ref{lem:subset} and Corollary~\ref{cor:T}.  Let 
   \begin{equation*}
\GABocd=\bigcup_{i=1}^q \GABbcdn g_i \GABbcdn
\end{equation*}
be a double coset decomposition of $\GABocd$ obtained by using the subgroup $G(\A,\bb,\B,\ddd,\cc)$. 
Now, as discussed in~\cite{Bremner2014}, either $$\left(\GABbcdn g_i \GABbcdn\right) \bigcap \Glpp=\emptyset$$
 or $$\GABbcdn g_i \GABbcdn \subseteq \Glpp.$$
Let $G_{\rm ext}$ be as in Step~\ref{step:Gext1} or Step~\ref{step:Gext2}. Then the set $\mathcal{T}^{L}_{{\rm Fix}_{G_{\rm ext}}}\neq \emptyset$
 if and only if the LP in Step~\ref{step:feasible} is feasible.
If $\mathcal{T}^{L}_{{\rm Fix}_{G_{\rm ext}}}\neq \emptyset$, then $G_{\rm ext} \leq \Glpp$ by Theorem~\ref{thm:thm}. Hence, $\GABbcdn$  can be updated with $G_{\rm ext}$.\qed
\end{pf}

\begin{algorithm}[t]
\centering
\caption{\label{meth:Glp} Computing  $\Glp$ of an LP $L$ of form~(\ref{eqn:geneqILP}) in standard form}
\begin{algorithmic}[1]
\State {\bf Input} $\A,\B,\bbb,\ddd,\cc$ from an  LP $L$ of form~(\ref{eqn:geneqILP}) in standard form.   
\State {\bf Compute}  the reduced row echelon form of $[\A\, |\, \bb]$ using Gaussian elimination;
\Comment{$\A\in \mathbb{R}^{m\times n}$.}\label{meth:Basicvar}
\State {\bf Substitute} the expressions obtained in Step~\ref{meth:Basicvar} for the basic variables in $\cc^{\top}\x$; \label{step:calculatex1}
\State {\bf Set}
$\hat{\cc}^{\top}\hat{\x}+a$ to be the resulting objective function from Step~\ref{step:calculatex1};\Comment{$ a \in \mathbb{R}$, $\hat{\x}\in \mathbb{R}^{n-m}$.}\label{step:chat}
\State {\bf Compute} $\Glpp$ and $ G(\A,\bb,\B,\ddd,\0)^{\rm LP}$ by using Method~\ref{meth:Glpp};
\State  {\bf Compute} $q$ and a set $\{g_1,\ldots,g_q\}$ so that  
$$
G(\A,\bb,\B,\ddd,\0)^{\rm LP}=\bigcup_{i=1}^q \Glpp g_i \Glpp;
$$\label{step:doublecoset}
\State {\bf Set} $\Glp:=\Glpp$;
\For{$i:= 1$ {\bf to} $q$ {\bf step} $1$}
\State {\bf Substitute} the expressions obtained in Step~\ref{meth:Basicvar} for the basic variables in $\cc^{\top}g_i(\x)$;
\label{step:calculatex2}
\State {\bf Set}
$\hat{\cc}_i^{\top}\hat{\x}+a_i$ to be the resulting objective function from Step~\ref{step:calculatex2};\Comment{$ a_i \in \mathbb{R}$, $\hat{\x}\in \mathbb{R}^{n-m}$.}\label{step:chati} 
\If{$\hat{\cc}^{\top}=\hat{\cc}_i^{\top}$ and $a=a_i$} \label{step:checkci}
\State {\bf Update} $\Glp:=\langle g_i,  \Glp \rangle$; 
\EndIf
\EndFor
\State {\bf Output} $\Glp$.\label{meth:GlpOut}
\end{algorithmic}
\end{algorithm}

 Method~\ref{meth:Glp} finds the $\Glp$ of an LP $L$ of form~(\ref{eqn:geneqILP}) in standard form by using the subgroup  $\Glpp$ of $\Glp$
  and the  symmetry group $G(\A,\bb,\B,\ddd,\0)^{\rm LP}$ of the feasible set of $L$. Both  $\Glpp$ and $G(\A,\bb,\B,\ddd,\0)^{\rm LP}$ 
  can be found by using Method~\ref{meth:Glpp}.
    Method~\ref{meth:Glp}
   requires computing a double coset decomposition of 
  $G(\A,\bb,\B,\ddd,\0)^{\rm LP}$ with respect to $\Glpp$. 
  This can be done by using {\tt GAP}~\cite{GAP2013}.
  The following lemma is used in proving the  theorem that establishes the viability of Method~\ref{meth:Glp}.
 \begin{lem}\label{lem:affine}
 Let $\A$ be an $m \times n$ matrix such that $m \leq n$ and 
 ${\rm rank}(\A)=m$. Let $\uu_i \in {\rm Null}(\A)$ for $i=1,\ldots,r$ be  linearly independent, where $r\leq n-m$.
Let $\ell$ be the set of indices of the basic variables 
in the reduced row echelon form of $[\A\, |\, \0]$.
 For $i=1,\ldots,r$, let $\uuuu_i\in 
 \mathbb{R}^{n-m}$ be obtained from $\uu_i$ by deleting its entries whose indices are in $\ell$.
 Then the vectors in $\{\uuuu_1,\ldots,\uuuu_{r}\}$ are linearly independent.
 \end{lem}
 \begin{pf} It suffices to show that the columns of $[\uuuu_1\,\uuuu_2\,\cdots\,\uuuu_r]$ are linearly independent, i.e., $${\rm rank}([\uuuu_1\,\uuuu_2\,\cdots\,\uuuu_r])=r.$$
 Since  $\uu_i \in {\rm Null}(\A)$, each entry  in $\uu_i$ whose index is  in $\ell$ is a linear combination of the entries whose indices are 
 in $\{1,\ldots,n-m\}\backslash \ell$.
 Then  $[\uu_1\, \uu_2\, \cdots\,\uu_r]$ is row equivalent to 
 $[\widetilde{\uu}_1\, \widetilde{\uu}_2\, \cdots\,\widetilde{\uu}_r]$, where $\widetilde{\uu}_i$ is obtained from $\uu_i$ by replacing each of its entries whose index is in $\ell$ with $0$. Now, since ${\rm rank}([\uu_1\, \uu_2\, \cdots\,\uu_r])={\rm rank}([\widetilde{\uu}_1\, \widetilde{\uu}_2\, \cdots\,\widetilde{\uu}_r])$ and 
 ${\rm rank}([\widetilde{\uu}_1\, \widetilde{\uu}_2\, \cdots\,\widetilde{\uu}_r])={\rm rank}([\uuuu_1\,\uuuu_2\,\cdots\,\uuuu_{r}])$, we get ${\rm rank}([\uu_1\, \uu_2\, \cdots\,\uu_r])=
 {\rm rank}([\uuuu_1\,\uuuu_2\,\cdots\,\uuuu_r])$.
 Hence,  $$r={\rm rank}([\uu_1\, \uu_2\, \cdots\,\uu_r])=
 {\rm rank}([\uuuu_1\,\uuuu_2\,\cdots\,\uuuu_r])$$
 as the vectors in $\{\uu_1,\ldots,\uu_r\}$ are linearly independent. \qed
 \end{pf}
The following theorem validates Method~\ref{meth:Glp}.
\begin{thm}
The output of Method~\ref{meth:Glp} is the $\Glp$ of $L$.
\end{thm} 
\begin{pf}
Let the matrix $\A$ be $m\times n$ with ${\rm rank}(\A)=m$.
Throughout the proof, for a feasible point $\vv$   of $L$, let $\vvvv$ be obtained from $\vv$ by deleting its entries whose indices are the same as those of the basic variables in the reduced row echelon form of $[\A\,|\,\bb]$, equivalently of $[\A\,|\,\0]$. 
Since the feasible set of $L$ is a full-dimensional 
polytope in an affine space of dimension $n-m$, the feasible set of 
$L$ contains affinely independent points $\vv_j\in \mathbb{R}^n$ for $j=0,\ldots,n-m.$ Then the vectors $\vv_j-\vv_0 \in {\rm Null}(\A)$
for $j=1,\ldots,n-m$  are linearly independent. 
 Hence,  the vectors in $\{\vvvv_1-\vvvv_0,\ldots, \vvvv_{n-m}-\vvvv_0\}$ are also linearly independent by Lemma~\ref{lem:affine}.

 Let $g_i$ be as in Step~\ref{step:doublecoset} of Method~\ref{meth:Glp}.  Let $\hat{\cc}^{\top}\hat{\x}+a$ and $\hat{\cc}_i^{\top}\hat{\x}+a_i$  be as in Step~\ref{step:chat} and 
 Step~\ref{step:chati} of Method~\ref{meth:Glp}.
First, we have $\Glp\leq 
G(\A,\bb,\B,\ddd,\0)^{\rm LP}$ as $G(\A,\bb,\B,\ddd,\0)^{\rm LP}$ is the  same as the symmetry group of the feasible set of $L$. Then by Definition~\ref{def:GAbBdc},
 $$ \Glpp \leq \Glp\leq 
G(\A,\bb,\B,\ddd,\0)^{\rm LP}. $$
Since either
$$ \Glpp g_i \Glpp \cap \Glp =\emptyset$$
or $$ \Glpp g_i \Glpp \subseteq \Glp, $$
it suffices to prove that $g_i \in \Glp$
if and only if $\hat{\cc}=\hat{\cc}_i$ and $a=a_i$. 

 Assume  $g_i\in\Glp$.
Then for  $j=0,\ldots,n-m$, we have $$\cc^{\top}\vv_j=\hat{\cc}^{\top}\vvvv_j+a, \quad \cc^{\top}g_i(\vv_j)=
\hat{\cc}_i^{\top}\vvvv_j+a_i,\quad \cc^{\top}\vv_j=\cc^{\top}g_i(\vv_j),$$ and we get
 $$\cc^{\top}\vv_j-\cc^{\top}\vv_0 = \cc^{\top}g_i(\vv_j)- \cc^{\top}g_i(\vv_0).$$ Hence, 
$$\hat{\cc}^{\top}\vvvv_j+a-\hat{\cc}^{\top}\vvvv_0-a = 
\hat{\cc}_i^{\top}\vvvv_j+a_i- \hat{\cc}_i^{\top}\vvvv_0-a_i,$$
and consequently,
$$\hat{\cc}^{\top}(\vvvv_j-\vvvv_0) = \hat{\cc}_i^{\top}(\vvvv_j-\vvvv_0).$$
Then,
\begin{equation}\label{eqn:cmct}
(\hat{\cc}-\hat{\cc}_i)^{\top}(\vvvv_j-\vvvv_0) = (\vvvv_j-\vvvv_0)^{\top}(\hat{\cc}-\hat{\cc}_i)=0
\end{equation}
for $j=1,\ldots,n-m$. Equations~(\ref{eqn:cmct}) imply 
 $$[(\vvvv_1-\vvvv_0)\,(\vvvv_2-\vvvv_0)\,\cdots\, (\vvvv_{n-m}-\vvvv_0)]^{\top}(\hat{\cc}-\hat{\cc}_i)=\0.$$
 Now, the linear independence of $\vvvv_j-\vvvv_0$ for $j=1,\ldots,n-m$ implies that the square matrix $$\left[(\vvvv_1-\vvvv_0)\,(\vvvv_2-\vvvv_0)\,\cdots\, (\vvvv_{n-m}-\vvvv_0)\right]^{\top}$$ is invertible. So, we conclude that $\hat{\cc}-\hat{\cc}_i=\0$ and $\hat{\cc}=\hat{\cc}_i$. Moreover, since 
 $$\hat{\cc}^{\top}\vvvv_0+a=\cc^{\top}\vv_0 = \cc^{\top}g_i(\vv_0)=\hat{\cc}_i^{\top}\vvvv_0+a_i,$$ we get $a=a_i$.
 
To prove the converse, assume $\hat{\cc}=\hat{\cc}_i$ and $a=a_i$. Since $g_i \in G(\A,\bb,\B,\ddd,\0)^{\rm LP}$, $g_i$ sends feasible points of $L$ to feasible points. 
Now, $\hat{\cc}=\hat{\cc}_i$ and $a=a_i$ implies that 
$$\cc^{\top}\vv=\hat{\cc}^{\top}\hat{\vv}+a=
\hat{\cc}_i^{\top}\hat{\vv}+a_i=\cc^{\top}g_i(\vv)$$ for each feasible point $\vv$ of $L$.
Hence, $g_i \in \Glp$ as $g_i$ preserves the feasibility and the objective function value of each feasible point.
\qed
\end{pf}
\begin{cor}
The symmetry group $\Glp$ of  an LP $L$ coincides with its formulation symmetry group  if the feasible set of  $L$ is non-empty,  
full dimensional, and $L$ has no redundant inequalities.
\end{cor}
\begin{pf}
Since the feasible set of $L$ is non-empty and full dimensional, there exists no equality constraint satisfied by each feasible point of $L$. Then, 
 $L$ has  no redundant inequality  constraints implies that $L$ is in standard form. WLOG assume that $L$ has the form of LP~(\ref{eqn:geneqILP}).
Then, $$G(\A,\bb,\B,\ddd,\cc)^{\rm LP}=\GABocd=G(\B,\ddd,\cc)=G(\A,\bb,\B,\ddd,\cc).$$ Let $L$
be the input to
Method~\ref{meth:Glp}. Then  $\hat{\cc}^{\top}\hat{\x}+a$ becomes $\cc^{\top}\x$ in
Step~\ref{step:chat},  $\hat{\cc}_i^{\top}\hat{\x}+a_i$ becomes 
$\cc^{\top}g_i(\x)=(g_i^{-1}(\cc))^{\top}\x$ in
Step~\ref{step:chati}, and the check in
Step~\ref{step:checkci} becomes $\cc^{\top}=(g_i^{-1}(\cc))^{\top}$.
So, the check in Step~\ref{step:checkci} requires that each new 
$g_i$ to be added to $G(\B,\ddd,\cc)$ must preserve 
the vector $\cc$. Then the output $\Glp$ must also preserve $\cc$. Hence, by Definition~\ref{def:GAbBdc}, $\Glp=G(\A,\bb,\B,\ddd,\cc)^{\rm LP}=G(\B,\ddd,\cc).$ \qed
\end{pf}
\section{Orthogonal arrays and their symmetries}
\label{sec:OA}
We first define orthogonal arrays (OAs). 
\begin{defn}\label{def:OA}
    An OA$(N, k, s, t)$ of strength $t \in \{0,\ldots,k\}$ is an $N \times k$ array of symbols from the set 
    $\{l_0,\ldots,l_{s-1}\}$   
 such that  each of the $s^t$  $t$-tuples  from 
 $\{l_0,\ldots,l_{s-1}\}^{t}$ appears $\lambda=N/s^{t}$ times in every $N\times t$ subarray. 
 \end{defn}
By Definition~\ref{def:OA}, every $N \times k$ array with symbols from a set $\{l_0,\ldots,l_{s-1}\}$ is an OA$(N,k,s,0)$ and vice versa. 
For  fixed $N$, $k$, $s$, and $t \in \{0,\ldots,k\}$, $\oan$ have many inherent symmetries, where each symmetry is a bijective map
from the set of all OA$(N,k,s,0)$ to  the set of all OA$(N,k,s,0)$ that preserves the $\oan$ property. In particular, each row permutation is a symmetry of $\oan$ for all $t\in \{0,\ldots,k\}$. We call each such  symmetry a {\it trivial symmetry} of  OAs.

A major source of non-trivial symmetries of $\oan$ 
 for  $t\in \{0,\ldots,k-1\}$ is the set of  isomorphism operations~\cite{Bulutoglu2016}. 
 (For $t=k$, it is easy to show that every isomorphism operation is a trivial symmetry.)
Next, we define isomorphism operations and the group of isomorphism operations that act  on $\oan$ for  $t\in \{0,\ldots,k\}$.
\begin{defn} Each of the $k!(s!)^k$ operations that involve permuting columns and the elements of $\{l_0,\ldots,l_{s-1}\}$ within each column of an $N \times k$ array with symbols from $\{l_0,\ldots,l_{s-1}\}$ is called an {\it isomorphism operation}. The set of all isomorphism operations forms a group called the {\it paratopism group}~\cite{Egan2016}. 
\end{defn}
We denote the paratopism group acting on $\oan$  with $G^{\rm iso}(k,s)$.
Two $\oan$s $\X$ and $\Y$ are {\it isomorphic} if 
$$
 \text{the set of rows of $\X$}= \text{the set of rows of $g(\Y)$} $$
  for some $g \in G^{\rm iso}(k,s)$~\cite{Stufken2007}. 
It is well known that $G^{\rm iso}(k,s) \cong S_s \wr S_k$~\cite{Egan2016}, where $S_s\wr S_k$ is the wreath product of the symmetric group of degree $s$ and the symmetric group of degree $k$. For a definition of the wreath product of groups, see~\cite{Rotman1994}. 
In~\cite{Bulutoglu2016}, OD-equivalence of $\oant$ for even $t$  was defined and used  to classify all non-isomorphic OA$(160,k,2,4)$ and OA$(176,k,2,4)$. To define  OD-equivalence of OAs, we  first need the concept of Hadamard equivalence from~\cite{McKay1979}.

\begin{defn}
Two $N\times k$ arrays $\Y_1$ and $\Y_2$ with symbols from $\{-1,1\}$ are {\it Hadamard equivalent} if $\Y_2$ can be obtained from $\Y_1$ by applying a sequence of signed permutations (permutations that may or may not be followed by sign changes) to the columns or rows
 of $\Y_1$.
\end{defn}

\begin{defn}
Two $N\times k$ arrays $\X_1$ and $\X_2$ with symbols from $\{-1,1\}$ are {\it OD-equivalent} if $[\1,\X_1]$ and $[\1,\X_2]$ are Hadamard equivalent.
\end{defn}
Clearly, two isomorphic $\oant$  with symbols from $\{-1,1\}$ are OD-equivalent. However, there exist OD-equivalent 
$\oant$ that are not isomorphic~\cite{Bulutoglu2016}. 
In what follows, we describe the operations other than the isomorphism operations that send an $\oant$ to one of its OD-equivalent copies.
Let $\Y$ be an $N \times k$ array with symbols 
from $\{-1,1\}$. For each $i \in \{1,\ldots,k \}$, define the column operation $R_i$ on $\Y$ to be

\begin{equation*}\label{eqn:M}
\Y=\begin{bmatrix}
\y_1& \cdots &  \y_i&\cdots&\y_k\\
 \end{bmatrix} \quad  \stackrel{R_i}{\longrightarrow}\quad
 \begin{bmatrix}
\y_1\odot\y_i& \cdots&\y_{i-1}\odot\y_i& \y_i&\y_{i+1}\odot\y_i&\cdots & \y_k\odot \y_i\\
 \end{bmatrix},
\end{equation*}
where 
$$\uu\odot\vv=\begin{bmatrix}
u_1v_1 \\
\vdots \\
u_n v_n\end{bmatrix}$$
for $\uu, \vv \in  \{-1,1\}^n$. 
Now, we have the following definition.
\begin{defn}\label{def:Ri}
Let the group generated by $R_1,\ldots,R_k$ and the elements of $G^{\rm iso}(k,2)$ be denoted by $G(k)^{\rm OD}$. Each element of $G(k)^{\rm OD}$ is called an {\it OD-equivalence  operation}.
\end{defn}
The proof of the next result is  a modification of the proof of Theorem 2 in~\cite{Arquette2016}. It also fills in  the details skipped in~\citep{Arquette2016}. 
We present Theorem 2 in~\cite{Arquette2016} as Theorem~\ref{thm:Arquette} in~Section~\ref{sec:OAapply}.
\begin{lem} \label{lem:finishOD}
Let $\Y=[\y_1\,\y_2\,\cdots\y_k] \in \{-1,1\}^{N \times k}$,  and $G(k)^{\rm OD}=\langle G^{\rm iso}(k,2),R_1,\ldots,R_k \rangle$ act on $\Y$ as in Definition~\ref{def:Ri}. Then  $G(k)^{\rm OD} \cong S_2^k\rtimes S_{k+1}$.
\end{lem}
\begin{pf}
Let 
$R= \langle R_1, \ldots, R_k \rangle \leq G(k)^{\rm OD}$  and $R_0=e$ be the identity element of $G(k)^{\rm OD}$. 
Let $\hat{S}_k$ be the group of all permutations that permute 
the columns of $\Y$.
Since   $R_iR_jR_i$ permutes $\y_i$ and $\y_j$, we have 
$\hat{S}_k \leq G(k)^{\rm OD}$ and $R_iR_jR_i\hat{S}_k=\hat{S}_k$. Then, $R_iR_jR_i\hat{S}_k=\hat{S}_k$ implies  $R_jR_i\hat{S}_k=R_i\hat{S}_k$ as 
$R_i^2=e$.
Since $R_jR_i\hat{S}_k=R_i\hat{S}_k$ for each distinct $i,j \in \{1,\ldots,k\}$,  there are $k+1$ left cosets of $\hat{S}_k$ in $R$. So each element  $x \in R$
can be written uniquely as $x=R_j\tau$ for some $\tau \in \hat{S}_k$ and $j \in \{0,\ldots,k\}$. Consequently,  $|R|=(k+1)!$. Then, $|R|<\infty$ implies that every element of $R$ can be written as a finite product of the $R_i$s.
Let 
\begin{equation}\label{eqn:trans}
\sigma_{ij}=R_iR_jR_i\in \hat{S}_k
\end{equation}
 for each distinct $i,j \in \{1,\ldots,k\}$. Then, for each distinct $i_1,\ldots, i_r \in \{1,\ldots,k\}$,
\begin{equation}\label{eqn:Rloop}
R_{i_1}R_{i_2}\cdots R_{i_{r-1}}R_{i_r}R_{i_1}=R_{i_1}R_{i_2}R_{i_1}R_{i_1}\cdots R_{i_1}R_{i_1}R_{i_{r-1}}R_{i_1}R_{i_1}R_{i_r}R_{i_1}=\sigma_{i_1i_2}\sigma_{i_1i_3}\cdots\sigma_{i_1i_r}
\end{equation}
\begin{equation}\label{eqn:Rnoloop}
R_{i_1}R_{i_2}\cdots R_{i_r}=R_{i_1}R_{i_r}R_{i_r}R_{i_2}R_{i_r}R_{i_r}\cdots R_{i_r}R_{i_r}R_{i_{r-1}}R_{i_r}=R_{i_1}R_{i_r}\sigma_{i_ri_2}\sigma_{i_ri_3}\cdots \sigma_{i_ri_{r-1}}
\end{equation}
by equation~(\ref{eqn:trans}) and $R_i^2=e$.
Now, since $R_{i_1}R_{i_r}=R_{i_r}\sigma_{i_1i_r}$,  equation~(\ref{eqn:Rnoloop})
becomes 
\begin{equation}\label{eqn:Rnoloopfin}
R_{i_1}R_{i_2}\cdots R_{i_r}=R_{i_r}\sigma_{i_ri_1}\sigma_{i_ri_2}\sigma_{i_ri_3}\cdots \sigma_{i_ri_{r-1}}.
\end{equation}
Then, given $x=R_{j_1}R_{j_2}\cdots R_{j_p}\in R$ for some not necessarily distinct $j_1,\ldots, j_p \in \{1,\ldots,k\}$, $x$ can be reduced to its unique form  $x=R_j\tau$ by first applying equation~(\ref{eqn:Rloop}) from right to left and then applying equation~(\ref{eqn:Rnoloopfin}) 
(if applicable) once equation~(\ref{eqn:Rloop})
can no longer be applied.

Given $x=R_{j_1}R_{j_2}\ldots R_{j_p} \in R$ for some not necessarily distinct $j_1,\ldots, j_p \in \{1,\ldots,k\}$, define $\psi:R\rightarrow S_{k+1}$ by 
$$\psi(x)=\psi(R_{j_1}R_{j_2}\cdots R_{j_p})=(j_1,k+1)(j_2,k+1)\cdots (j_p,k+1).$$ We now prove that $\psi$ is an isomorphism. First, assuming that $\psi$ is well-defined, it is clear that $\psi$ is a homomorphism. Second, $S_{k+1}$ is generated by the elements of $\{(1,k+1),(2,k+1),\ldots, (k,k+1)\}$ as  each transposition $(i,j)$ satisfies
\begin{equation}\label{eqn:trans2}
(i,j)=(i,k+1)(j,k+1)(i,k+1).
\end{equation}
Thus, 
 $\psi$ is onto $S_{k+1}$. Assuming that $\psi$ is a well-defined map, injectivity of $\psi$ follows from the facts that $|S_{k+1}|=|R|$
 and  $\psi$ is onto $S_{k+1}$. 
Hence, it suffices to show that $\psi$ is a well-defined map.
Let $R'_i=(i,k+1)$ for $i \in \{1,\ldots, k\}$ and $\sigma'_{ij}=(i,j)\in S_k$ for each distinct $i,j \in \{1,\ldots, k\}$. Let $R'_0=e'$ be the identity permutation in $S_{k+1}$.
Then equation~(\ref{eqn:trans2})
and $R'^2_i=e'$ imply
\begin{align}
R'_{i_1}R'_{i_2}\cdots R'_{i_{r-1}}R'_{i_r}R'_{i_1}=&\ 
\sigma'_{i_1i_2}
\sigma'_{i_1i_3}\cdots\sigma'_{i_1i_r} \label{eqn:Rloop2} \\
R'_{i_1}R'_{i_2}\cdots R'_{i_r}=&\ R'_{i_r}\sigma'_{i_ri_1}\sigma'_{i_ri_2}\sigma'_{i_ri_3}\cdots \sigma'_{i_ri_{r-1}}\label{eqn:Rnoloopfin2}
\end{align} 
the same way equation~(\ref{eqn:trans}) and $R_i^2=e$ imply equations~(\ref{eqn:Rloop}) and~(\ref{eqn:Rnoloopfin}).
Let $x=R_{j_1}R_{j_2}\cdots R_{j_{p_1}} \in R$ and $y=R_{j'_1}R_{j'_2}\cdots R_{j'_{p_2}}\in R$ be such that $x=y$. To finish the proof, we need to show that 
$\psi(x)=\psi(y)$. 
The equality $x=y$ implies that $x=y=R_{j} \tau$ for some $j \in \{0,\ldots,k\}$ and $\tau \in \hat{S}_k$. Moreover,  $x=y=R_{j} \tau$ can be obtained  by first applying equation~(\ref{eqn:Rloop}) 
from right to left and then applying equation~(\ref{eqn:Rnoloopfin}) (if applicable) to $x=R_{j_1}R_{j_2}\cdots R_{j_{p_1}}$ and $y=R_{j'_1}R_{j'_2}\cdots R_{j'_{p_2}}$. Then $\psi(x)=\psi(y)=R'_j\tau'$, for some $\tau'\in S_k$ can be obtained  by first applying equation~(\ref{eqn:Rloop2}) 
from right to left and then applying equation~(\ref{eqn:Rnoloopfin2}) (if applicable) to $\psi(x)=R'_{j_1}R'_{j_2}\cdots R'_{j_{p_1}}$ and $\psi(y)=R'_{j'_1}R'_{j'_2}\cdots R'_{j'_{p_2}}$.
Therefore, we conclude that $R \cong S_{k+1}$.

 Let $\phi \in Z_2^k$ be a multiplication of some subset of columns of $\Y$ by $-1$. Then $R_i \phi R_i=R_i^{-1} \phi R_i=\phi'$, where  
$\phi' \in Z_2^k$. This implies that $Z_{2}^k\trianglelefteq G(k)^{\rm OD}$ as $Z_2^k \trianglelefteq Z_2\wr \hat{S}_k \leq G(k)^{\rm OD}$. Then, $Z_{2}^k\trianglelefteq G(k)^{\rm OD}$, $Z_2^k \cap R=\{e\}$, and  $R\leq G(k)^{\rm OD}$ imply 
$Z_2^k\rtimes R \leq G(k)^{\rm OD}$.
Now, $G(k)^{\rm OD}=\langle R,Z_2^k \rangle$ and $Z_2^k\trianglelefteq G(k)^{\rm OD}$. So, for each $g \in  G(k)^{\rm OD}$,  $g=rw$ for some $r \in R$ and 
$w\in Z_2^k$. Consequently, $|G(k)^{\rm OD}|\leq |R||Z_2^k|=(k+1)!2^k$. Hence, we get $G(k)^{\rm OD}=Z_2^k\rtimes R \cong S_2^k\rtimes S_{k+1}.$\qed
\end{pf}

For even $t$, a major source of non-trivial symmetries of $\oant$ that are not
isomorphism operations are the OD-equivalence operations that are not in $G^{\rm iso}(k,2)$~\cite{Bulutoglu2016}. The  following theorem shows that the OD-equivalence operations are indeed symmetries of $\oant$ when $t$ is even.
\begin{thm}\label{thm:OD}
Let  $\Y$ be an $\oant$  with symbols from $\{-1,1\}$ and strength $t \geq 1$.  Then $\X$ is OD-equivalent to $\Y$ if and only if  there exists an OD-equivalence operation $g$ such that $$\text{the set of rows of $\X$}=\text{the set of rows of $g(\Y)$}.$$ Moreover, if $\X$ is OD-equivalent to $\Y$, then $\X$ is
an OA$(N,k,2,2\lfloor \frac{t}{2} \rfloor)$.  
\end{thm}
We first show that two OAs each with strength at least $1$ are OD-equivalent if and only if the set of rows of one can be obtained from that of the other by applying an element of $G(k)^{\rm OD}$ to each row in the set. 
\begin{lem}\label{lem:PID}
Let $\X$, $\Y$ be  OA$(N,k,2,t_1)$,  OA$(N,k,2,t_2)$ with symbols from $\{-1,1\}$ and strengths $t_1\geq 1$, $t_2 \geq 1$. Then  $\X$ and $\Y$ are OD-equivalent
if and only if  $$
 \text{the set of rows of $\X$}= \text{the set of rows of $g(\Y)$} $$
  for some $g \in G(k)^{\rm OD}$.  
\end{lem}
\begin{pf}
By the definition of OD-equivalence, $\X$ and $\Y$ are OD-equivalent if and  only if  
 \begin{equation}\label{eqn:ODfirst}
\ppi\D[\1\,\X]=[\1\,\Y]\D_1\ppi_1
\end{equation} for some  permutation matrices  
$\ppi$ and 
$\ppi_1$ and diagonal matrices $\D$ and $\D_1$  whose diagonal entries are in $\{-1,1\}$. 
Now, $\Y$ is an OA with strength at least $1$; consequently, $\ppi\D[\1\,\X]=[\1\,\Y]\D_1\ppi_1$ must have $\pm \1$ as a column exactly once. Since $\X$ is also an OA with strength at least $1$, we must have either $\D=\pm\I$ or $\D=\pm \mbox{diag}(\x_{i})$, where   $\x_i$ is the $i$th column of $\X$ for some $i$. 
When $\D=\pm\I$, equation~(\ref{eqn:ODfirst}) holds if and only if  $\X$ is isomorphic to $\Y$. This is true if and only if $\X$ can be obtained from $\Y$ by applying an element $g$ of $G^{\rm iso}(k,2)\leq G(k)^{\rm OD}$.
When $\D=\pm \mbox{diag}(\x_{i})$ for some $i$,
 then equation~(\ref{eqn:ODfirst}) holds if and only if
 \begin{equation}\label{eqn:SetM}
 \{\ppi\D[\1\,\X]\} \backslash \{\1, -\1 \}=\{[\1\,\Y]\D_1\ppi_1\} \backslash \{ \1, -\1\},
\end{equation} 
  where $\{\D[\1\,\X]\}\backslash \{\1,-\1\}=\{\pm\x_{i}\odot\x_{1}, \ldots, \pm\x_{i}\odot\x_{i-1}, \pm\x_{i}, \pm\x_{i}\odot\x_{i+1}, \ldots,  \pm\x_{i}\odot \x_{k}\}$ and
$\{\M\}$ is the set of columns of a matrix $\M$. 
 Equation~(\ref{eqn:SetM}) holds if and only if $$
 \text{the set of rows of $h_1(\X)$}= \text{the set of rows of $h_2(\Y)$} $$
  for some $h_1,h_2 \in G(k)^{\rm OD}$.
 The result now follows from 
$$\text{the set of rows of $\X$}= \text{the set of rows of $h_1^{-1}h_2(\Y)$} $$ by taking $g=h_1^{-1}h_2$.
 \qed
\end{pf}

Next, we prove Theorem~\ref{thm:OD},
 but first  we need 
 the concept of $J$-characteristics and the subsequent lemma from~\cite{Stufken2007}. 
\begin{defn}
Let $\Y=[y_{ij}]$ be an $N \times k$ array with symbols from 
$\{-1,1\}$. Let $r \in \{1,\ldots,k\}$ and $\ell=\{i_1,\ldots,i_r\}\subseteq \{1,\ldots,k\}$.   Then the integers
\begin{equation*}
J_r(\ell)(\Y):=\sum_{i=1}^N{\prod_{j \in \ell}{y_{ij}}}
\end{equation*}
are called the {\it $J$-characteristics} of $\Y$. (For $r=0$, $J_0(\emptyset)(\Y):=N$.)
\end{defn}
\begin{lem}[Stufken and Tang~\cite{Stufken2007}]\label{lem:J}
An $N\times k$ array $\Y$ with symbols from  $\{-1,1\}$ is an $\oant$ if and only if  $J_r(\ell)(\Y)=0$ for all $\ell \subseteq \{1,\ldots,k\}$ such that $|\ell|=r$ and $r\in \{1,\ldots,t\}$.
 \end{lem}
 \begin{lem}\label{lem:Jr}
Let   $\ell \subseteq \{1,\ldots,k\}$ be such that $|\ell|=r>0$. 
Let $g$ be an OD-equivalence operation and $g(\Y)$ be the array obtained after $g$ is applied to $\Y$.
 Then $$J_r(\ell)(g(\Y))=\pm J_{r'}(\ell')(\Y)$$ for some $\ell' \subseteq \{1,\ldots,k\}$, where 
 \begin{eqnarray}\label{eqn:spann}
|\ell'|=r'=\left\{
\begin{array}{rl}
\text{$r$ or $r+1$} & \quad \text{if $r$ is odd,}\\
\text{$r$ or $r-1$} &\quad  \text{otherwise.}
\end{array}\right.
\end{eqnarray}
 \end{lem}
\begin{pf}
Let $i\in\{1,\ldots,k\}$ and $R_i$ be as in Definition~\ref{def:Ri}. 
Then, 
\begin{equation}\label{eqn:Ri}
J_r(\ell)(R_i(\Y))= \begin{cases}
J_r(\ell)(\Y)\quad \quad \quad\ \quad \quad \hspace{.05cm} \text{if $r$ is even and $i \notin \ell$,}\\
J_{r-1}(\ell\backslash\{i\})(\Y)\quad \quad \, \, \, 
\hspace{.05cm} \text{if $r$ is even and $i \in \ell$,}\\
J_{r+1}(\ell\cup \{i\})(\Y)\quad \quad \text{if $r$ is odd and $i \notin \ell$,}\\
J_r(\ell)(\Y) \quad \quad \quad \quad \quad \hspace{-.03cm} \, \, \, \text{if $r$ is odd and $i \in \ell$.}
\end{cases}
\end{equation}
Let $R=\langle R_1,\ldots,R_k \rangle$ and $Z_2^k$ be the group of all possible sign switches of columns of $\Y$.
 Then by the proof of Lemma~\ref{lem:finishOD}, $g=g_1g_2$, where $g_1 \in R$   and $g_2 \in Z_2^k$. 
Then by equation~(\ref{eqn:Ri}),  
\begin{equation*}\label{eqn:isomorphism}
J_r(\ell)(g(\Y))=J_r(\ell)(g_1(g_2(\Y)))=J_{r'}(\ell')(g_2(\Y))
\end{equation*}
for some $\ell' \subseteq \{1,\ldots,k\}$ and $r'=|\ell'|$  as in equation~(\ref{eqn:spann}).
 Now, $g_2(\Y)$ is obtained from $\Y$ by multiplying some columns of $\Y$ by $-1$. Hence,
 $$J_r(\ell)(g(\Y))=J_{r'}(\ell')(g_2(\Y))=\pm J_{r'}(\ell')(\Y).$$
\qed
\end{pf}
Finally, observe that Theorem~\ref{thm:OD} follows from Lemmas~\ref{lem:PID}, \ref{lem:J}, and~\ref{lem:Jr}.
\section{The LP relaxation symmetry group of an $\oan$ defining ILP}\label{sec:OAapply}
First, for each $N,k,s,t$ combination, we describe the $\oan$ defining ILP in~\cite{Bulutoglu2008} whose   feasible set  contains a set of all non-isomorphic (non-OD-equivalent if $s=2$) $\oan$.
Define the {\em frequency vector} of an $\oan$ to be $\x:=(x_1,\ldots,x_{s^k})$ whose $\left(\sum_{j=1}^k i_j s^{k-j}+1\right)$th entry  is the number of times the  symbol combination $(l_{i_1},\ldots,l_{i_k}) \in \{l_0,\ldots,l_{s-1}\}^k$ appears in an $\oan$. 
  Then $\x$ must be a feasible point of  the  ILP 
\begin{align}
 &\quad \quad \quad \quad   \quad \min  \quad \1^{\top}\x \nonumber\\
 \text {s.t.}\hspace{0cm} & \sum_{\{i_{1},\ldots,i_{k}\} \backslash \{i_{j_1}, \ldots,i_{j_t}\} \in \{0,\ldots,{s-1}\}^{k-t} } {x_{[i_1 s^{k-1}+\cdots+i_k s^{k-k}+1]}}= \frac{N}{s^{t}} \label{ilp:BF} \\
 &\text{for each}\ \{j_1, \ldots,j_t\}\subseteq \{1,\ldots,k\}\ \text{and}\ (i_{j_1}, \ldots, i_{j_t}) \in \{0,\ldots,s-1\}^{t}, \nonumber\\
 &   0 \leq x_i \leq  p_{\max}, \quad x_i \in 
 \mathbb{Z}, \quad \text{for} \  i \in \{1,\ldots, s^k\} \nonumber
\end{align}
with a large formulation symmetry group, 
 where 
 $\1$ is the vector of all ones and $p_{\max}\leq N/s^{t}$ is a positive integer computed as in~\cite{Bulutoglu2008}. All feasible points (solutions) of ILP~(\ref{ilp:BF}) are optimal as each $\oan$ must have $N$ rows. So, the objective function 
 $\1^{\top}\x$ was introduced to formulate a constraint satisfaction problem for OAs as an ILP.
   
 Let $G(k,s,t)$ be the formulation symmetry group of  ILP~(\ref{ilp:BF}). 
In~\cite{Geyer2014}, it is shown that $G(k,s,t) \cong S_s\wr S_k$ for $1 \leq t\leq k-1$, and each element of $G(k,s,t)$ sends the frequency vector of an $\oan$ to that of one of its isomorphic copies. Hence, for $1 \leq t\leq k-1$, $G(k,s,t)= G^{\rm iso}(k,s)$ as $G^{\rm iso}(k,s)$'s action on the frequency vector of an $\oan$ is identical to that of $G(k,s,t)$. 
In~\cite{Bulutoglu2008}, all $\oan$ for many $N,k,s,t$ combinations were enumerated up to isomorphism   by finding a set of all non-isomorphic solutions to ILP~(\ref{ilp:BF}) under the action of $G(k,s,t)$. 

  Let $ G(k,s,t)^{\rm LP}$ be the LP relaxation symmetry group  
 of ILP~(\ref{ilp:BF}).
Since the formulation symmetry group of an LP is a subgroup of the LP relaxation symmetry group, we have the following result.
  \begin{lem}\label{lem:wr}
The LP relaxation symmetry group $ G(k,s,t)^{\rm LP}$ contains $G(k,s,t)=G^{\rm iso}(k,s) \cong S_s\wr S_k$, and hence $ |G(k,s,t)^{\rm LP}| \geq |S_s\wr S_k|=k!(s!)^k$.
\end{lem}
Having linearly dependent constraints in an ILP formulation
slows down a B\&B algorithm as LP relaxations in the B\&B search tree take longer to solve.  The  ILP formulation 
\begin{align}
 &\quad \quad \quad \quad   \quad \min  \quad 0 \nonumber\\
 \text {s.t.}\hspace{0cm} & \sum_{\{i_{1},\ldots,i_{k}\} \backslash \{i_{j_1}, \ldots,i_{j_q}\} \in \{0,\ldots,s-1\}^{k-q} } {x_{[i_1 s^{k-1}+\cdots+i_k s^{k-k}+1]}}=
  \frac{N}{s^q}\nonumber\\ &\text{for each}\    q \in \{0,\ldots, t\},\  \{j_1, \ldots,j_q\}\subseteq \{1,\ldots,k\},\label{ilp:BFimproved}\\& \text{and}\  (i_{j_1}, \ldots, i_{j_q}) \in \{0,\ldots,s-2\}^q, \nonumber \\
 & 1 \leq x_1, \quad 0 \leq x_i \leq  p_{\max}, \quad x_i \in 
 \mathbb{Z}, \quad \text{for} \  i \in \{1,\ldots, s^k\} \nonumber
\end{align}
from~\cite{Bulutoglu2016} improves the
 ILP~(\ref{ilp:BF}) formulation by replacing its set of linearly dependent equality constraints with a row equivalent, yet linearly independent set of equalities. (We are going to prove that ILP~(\ref{ilp:BFimproved}) has a linearly independent set of equality constraints.)
 ILP~(\ref{ilp:BFimproved}) has a smaller formulation symmetry group than that of ILP~(\ref{ilp:BF}), and yet the following remark holds.
 \begin{rem}\label{rem:sol}
 The set of  all lexicographically minimum  solutions of ILP~(\ref{ilp:BFimproved}) under the action of $G^{\rm iso}(k,s)=G(k,s,t)$ or $G(k,s,t)^{\rm LP}$ is the same as that of ILP~(\ref{ilp:BF}). 
\end{rem} 
 By Lemma~\ref{lem:wr}, we have $G^{\rm iso}(k,s)=G(k,s,t) \leq 
G(k,s,t)^{\rm LP}$. Hence, it suffices to justify Remark~\ref{rem:sol} for 
$G^{\rm iso}(k,s)$. Since  $G^{\rm iso}(k,s)$ acts transitively on the indices of the variables  of ILP~(\ref{ilp:BF}), all lexicographically minimum solutions of ILP~(\ref{ilp:BF}) under $G^{\rm iso}(k,s)$ satisfy $1 \leq x_1$.
Then  Remark~\ref{rem:sol} follows as the equality constraints of ILP~(\ref{ilp:BF}) and ILP~(\ref{ilp:BFimproved}) are row equivalent and ILP~(\ref{ilp:BF}) and 
ILP~(\ref{ilp:BFimproved}) have the same inequalities if $1 \leq x_1$ is deleted from ILP~(\ref{ilp:BFimproved}).
 Hence, finding a set of all non-isomorphic solutions to ILP~(\ref{ilp:BFimproved})  under the action of $G(k,s,t)$ is equivalent to classifying  $\oan$ up to isomorphism. 
 
 Let $\A(k,s,t)$ and $\A'(k,s,t)$ be the equality constraint matrices of ILPs~(\ref{ilp:BF}) and~(\ref{ilp:BFimproved}).
 Next, we are going  to establish a complete characterization of $G(k,s,t)^{\rm LP}$ that only involves the projection matrix $\PP_{\A(k,s,t)^{\top}}=\PP_{\A'(k,s,t)^{\top}}$ by proving the following theorem.
\begin{thm}\label{thm:Glp}
 Let $p_{\max}\neq \lambda/(s^{k-t})$ in ILP~(\ref{ilp:BF}). Then  $G(k,s,t)^{\rm LP}=G_{\PP_{\A'(k,s,t)^{\top}}}$, where $G_{\PP_{\A'(k,s,t)^{\top}}}$
 is the automorphism group of $\PP_{\A'(k,s,t)^{\top}}$.
\end{thm} 
We are also going to show  
\begin{equation}\label{eqn:Arquette}
\left|G(k,s,t)^{\rm LP}\right| \geq
\begin{cases}
\left(k+1\right)!2^k &  \mbox{if }  s=2 \mbox{ and $t$ is even, } \\
 k!\left(s!\right)^k & \text{otherwise}
\end{cases}
\end{equation}  
 for $1 \leq t \leq k-1$   and  prove the following theorem.
 \begin{thm}\label{thm:use}
Let  $1\leq t\leq k-1$ and $G(k,s,t)^{\rm LP}$ satisfy inequality~(\ref{eqn:Arquette}) as an equality. Let $G(k,s,t)^{\rm LP}$ be used within B\&B with isomorphism pruning  to find a set of all non-isomorphic solutions $\mathcal{F}$ of ILP~(\ref{ilp:BFimproved}). Then 
\begin{equation*}\label{eqn:Arquetteoa}
\mathcal{F} =
\begin{cases}
\text{a   set of all non-OD-equivalent } OA$(N,k,$s$,t)$&  \mbox{if   $t$ is even and $s=2$, } \\
 \text{a set of all non-isomorphic } $\oan$ & \text{otherwise.}
\end{cases}
\end{equation*}
\end{thm}
In~\cite{Arquette2016},  inequality~(\ref{eqn:Arquette}) was proven to be  an equality, i.e., the hypothesis of Theorem~\ref{thm:use} was proven when $s=2$, $t=1$ and when $s=t=2$, $k\geq 4$ by proving the following two theorems.
\begin{thm}[Arquette and Bulutoglu~\cite{Arquette2016}]\label{thm:Arquette1}
$G(k,2,1)^{\rm LP} \cong S_2^k \rtimes S_{k}$.
\end{thm}
\begin{thm}[Arquette and Bulutoglu~\cite{Arquette2016}]\label{thm:Arquette}
 For    $k \geq 4$, 
$G(k,2,2)^{\rm LP} \cong S_2^k \rtimes S_{k+1}$.
\end{thm}
Once a set of all non-OD-equivalent $\oant$ is found, the method in~\cite{Bulutoglu2016} for extracting a set of all non-isomorphic $\oant$ can be used to find a set of all non-isomorphic $\oant$. The time it takes for such an extraction 
 was observed to be insignificant compared to the time it takes to find a set of all non-OD-equivalent OA$(N,k,2,t)$~\cite{Bulutoglu2016}. 
 
Theorem~1 in~\cite{Geyer2014} modifies ILP~(\ref{ilp:BF})   to an ILP  without equalities  by deleting a set of basic variables after Gaussian elimination. Such a modification is only possible because at the end of Gaussian
elimination the coefficient of each basic variable  is one and the
coefficient of each free variable is an integer. The deleted variables are in fact a set of slack variables of the resulting ILP's LP relaxation. For the
 computational experiments in~\cite{Geyer2014}, the objective function of this ILP was taken to be the zero function. Let $G(k,s,t)^{\rm LP \leq}$ be the LP relaxation symmetry group of the resulting ILP once the basic variables are deleted. In~\cite{Geyer2014},  it is shown that $G(k,s,t)^{\rm LP \leq}\cong S_{s-1} \wr S_{k}$. Hence, $$\left|G(k,s,t)^{\rm LP\leq}\right|=k!((s-1)!)^k<\left|G(k,s,t)\right|=k!(s!)^k.$$ On the other hand, since ILP~(\ref{ilp:BF}) has equality constraints, it is not clear whether $G(k,s,t)=G(k,s,t)^{\rm LP}$.  
 
 Next, we describe yet another $\oant$ defining ILP formulation  developed in~\cite{Arquette2016}. We will add redundant equalities to this formulation to show that 
 $G(k,2,t)^{\rm LP}$ contains  $G(k)^{\rm OD}$ when $t$ is even.
 This ILP has the same variables, the same inequalities (excluding the inequality $1\leq x_1$ in ILP~(\ref{ilp:BFimproved})) as in  ILPs~(\ref{ilp:BF}) and~(\ref{ilp:BFimproved}). Moreover its  equality constraints are row equivalent to those of ILPs~(\ref{ilp:BF}) and~(\ref{ilp:BFimproved}). Hence, its LP relaxation's feasible set is the same as that of the LP relaxation of  ILP~(\ref{ilp:BF}).
 
  Let the transpose of row vectors of $\ZZ=[\z_1\, \z_2\, \cdots\, \z_k]$ be all $2^k$  vectors in $\{-1,1\}^k$. For   $i_1< \cdots < i_r \in \{1, \ldots, k\}$ with $r\geq2$, let  $\mathbf{z}_{i_1, \ldots, i_r}$  be the  $r$-way Hadamard product $\mathbf{z}_{i_1} \odot \cdots \odot \mathbf{z}_{i_r}$, where the $p$th entry of $\mathbf{z}_{i_1} \odot \cdots \odot \mathbf{z}_{i_r}$ is
the product of the entries on the $p$th row of the matrix   $[\mathbf{z}_{i_1}\, \mathbf{z}_{i_2}\, \cdots\,  \mathbf{z}_{i_r}]$.
Let $x_p$ for $p=1,\ldots,2^k$ be the number of times the $p$th row of $\ZZ$ appears in an $\oant$.
Now, by Lemma~\ref{lem:J}, the ILP

\begin{eqnarray}\label{eqn:Mf}
\min \ \1^{\top} \mathbf{x}\ \nonumber\\
\text{s.t.} \ \ \ \ \quad \1^{\top} \mathbf{x} = N, \nonumber\\ 
\begin{bmatrix}[r] 
\M \\ 
-\M
\end{bmatrix}\x  = \mathbf{\0},\\
 \mathbf{x} \geq \0,   \ \mathbf{x}\in \mathbb{Z}^{2^k} \nonumber
\end{eqnarray}
is an $\oant$ defining ILP formulation,
 where  $\mathbf{M}$  is the  $\sum\limits_{i=1}^{t} \binom{k}{i} \times 2^k$  matrix
\begin{equation*}\label{eqn:M}
\mathbf{M} = \begin{bmatrix}
\mathbf{z}_1^{\top} \\
\vdots \\
\mathbf{z}_k^{\top} \\
\mathbf{z}_{1, 2}^{\top} \\
\vdots \\
\mathbf{z}_{{k - t + 1}, \dots, k}^{\top} \end{bmatrix}.
\end{equation*}
The constraints 	$\1^{\top} \mathbf{x} = N$ and $\mathbf{M} \mathbf{x} = \mathbf{\0}$ ensure that the sought after OAs have $N$ rows and all of their $J_r(\ell)=0$
for all $\ell\subseteq\{1,\ldots,k\}$ such that $1 \leq |\ell|=r\leq t$. We added the redundant equalities $-\mathbf{M} \mathbf{x} = \0$ so that we can prove the following lemma. 
\begin{lem}\label{lem:contains}
Let $1 \leq t \leq k-1$ and  $\x$ be the frequency vector of an $\oant$ as in ILP~(\ref{eqn:Mf}).  If $t$ is even (odd),  then the formulation symmetry group of ILP~(\ref{eqn:Mf}) contains 
 $G(k)^{\rm OD}$ $(\Giso)$, where  $g(\x)$ is the frequency vector of an OD-equivalent (isomorphic) $\oant$ for each element $g$ of $G(k)^{\rm OD}$ $(\Giso)$.
\end{lem}	
\begin{pf}
The group $G(k)^{\rm OD}$ ($\Giso$) acts on $N \times k$ arrays with symbols from $\{-1,1\}$. This action induces an action on the frequency vectors of such arrays. Identify the action of $ G(k)^{\rm OD}$ ($\Giso$) on $\x$ by $g(x_{p_1})=x_{p_2}$ if and only if $$g((z_{p_11},\ldots, z_{p_1k}))=(z_{p_21},\ldots, z_{p_2k}),$$ where $p_1,p_2 \in\{1,\ldots,2^k\}$ and $g$ is an  OD-equivalence (isomorphism) operation on the columns of the row vector $(z_{p_11},\ldots, z_{p_1k})$. Hence, $g$ sends $\ZZ$ to one of its OD-equivalent
(isomorphic) copies. Here, $\x$ is indexed by the rows of $\ZZ$.
Since $g$ 
is an invertible map from $\{-1,1\}^k$ to $\{-1,1\}^k$ and
$$\text{the set of columns of}\,\, \ZZ^{\top}=\{-1,1\}^k,$$
$g$ is a permutation of the rows of $\ZZ$. Consequently, $g$ is  a permutation of columns of $\M$.  
Now, since $g$ is an OD-equivalence (isomorphism) operation on the columns of $\ZZ$, $g$'s action on the columns of $\M$ sends the row $\z_{i_1,\ldots, i_r}^{\top}$ to a row of the form 
$\pm \z_{i'_1,\ldots, i'_{r'}}^{\top}$ ($\pm \z_{i''_1,\ldots, i''_{r}}^{\top}$), where  $\{i_1,\ldots, i_r\},
\{i'_1,\ldots, i'_{r'}\}\subseteq \{1,\ldots,k\}$ ($\{i_1,\ldots, i_r\},
\{i''_1,\ldots, i''_{r}\}\subseteq \{1,\ldots,k\}$), and
\begin{eqnarray*}\label{eqn:rp}
r'=\left\{
\begin{array}{rl}
\text{$r$ or $r+1$} & \quad \text{if $r$ is odd,}\\
\text{$r$ or $r-1$} &\quad  \text{otherwise.}
\end{array}\right.
\end{eqnarray*}
Hence, when $t$ is even (odd), we get
\begin{eqnarray}\label{eqn:MfPi}
\begin{bmatrix}[r]
 \M \\ 
-\M
\end{bmatrix}g (\x) = 
\ppi\begin{bmatrix}[r]
\M \\
-\M
 \end{bmatrix} \x =\ppi\mathbf{\0}=\0,\nonumber 
\end{eqnarray}
 where $\ppi$ is a $\sum_{i=1}^{t} 2{k \choose i} \times  \sum_{i=1}^{t} 2{k \choose i} $ permutation matrix.  
 This proves that for even (odd) $t$ the formulation symmetry group of ILP~(\ref{eqn:Mf})  contains $G(k)^{\rm OD}$ 
 ($\Giso$)  as $g$ maps both the objective function $\1^{\top}\x$ and the constraint $\1^{\top}\x=N$ to themselves and permutes the constraints  $\x\geq\0$   among each other. 
 Finally, by Lemma~\ref{lem:PID} applied to the $\oant$ whose frequency vector is $\x$ for $t \in \{1,\ldots,k\}$ 
  (as the action of each element of $\Giso$ sends the frequency vector of an $\oant$ to that of one of its isomorphic copies), $g(\x)$ is the frequency vector of an OD-equivalent (isomorphic) $\oant$.
 \qed
\end{pf}
\begin{lem}\label{lem:containsnot}
For odd $t$ and $1\leq t\leq k-1$, $G(k,2,t)^{\rm LP}$ contains none of the $R_i$ in Definition~\ref{def:Ri}. 
\end{lem}	
\begin{pf}
As in the proof of Lemma~\ref{lem:contains}, identify the action of $G(k)^{\rm OD}$ on the frequency vectors $\x$ of $N\times k$ arrays with symbols from $\{-1,1\}$. 
Pick $1\leq i_1< \cdots <i_t \leq k $ and $i \in \{1,\ldots,k\} \backslash \{i_1,\ldots, i_t\}$.
Take $g$ in  Lemma~\ref{lem:contains} to be $R_i$  and observe that $R_i$ is a linear transformation from $\mathbb{R}^{2^k}$ to $\mathbb{R}^{2^k}$ as $R_i$ is a permutation of the coordinates of $\x$. Let $\R_i$ be the matrix of $R_i$ with respect to the standard basis.
Then $\M \R_i$ has the $(t+1)$-way Hadamard product $(\z_i\odot(\z_{i_1}\odot\cdots\odot\z_{i_t}))^{\top}$ as one of its rows, and this row is orthogonal to all the rows of $\M$.
Hence, the action of $R_i$ on the constraints of ILP~(\ref{eqn:Mf})
produces an equality constraint that is not a linear combination of the original constraints, so $R_i \notin G(k,2,t)^{\rm LP}$.\qed
\end{pf} 
The following theorem is in part stated, but not proven in~\cite{Arquette2016}.
\begin{lem}\label{lem:contains2}
Let $1\leq t\leq k-1$. Then $G(k,2,t)^{\rm LP}$ contains  $G(k)^{\rm OD}\cong S_2^k \rtimes S_{k+1}$, where each element of $G(k)^{\rm OD}$ sends the frequency vector of an $\oant$ to that of an OD-equivalent $\oant$ if and only if $t$ is even. Hence, for even $t$, $|G(k,2,t)^{\rm LP}|\geq |S_2^k \rtimes S_{k+1}|=(k+1)!2^k$.
\end{lem}
\begin{pf} 
The equality constraints of ILP~(\ref{eqn:Mf}) can be obtained as linear combinations of those of ILP~(\ref{ilp:BF}). Rows of $\M$ and $\1^{\top}$ form a mutually orthogonal set of vectors
 of size $\sum_{i=0}^{t} {k \choose i}$ in $\mathbb {R}^{2^k}$.
 Hence, the rank of the equality constraint matrix of ILP~(\ref{eqn:Mf}) is
  $\sum_{i=0}^{t} {k \choose i}$.
 By Lemma 1 in~\cite{Rosenberg1995}, $\sum_{i=0}^{t} {k \choose i}$ is also the rank 
 of the equality constraint matrix of ILP~(\ref{ilp:BF}).
 Now, this implies that the equality constraints of  ILPs~(\ref{ilp:BF}) and~(\ref{eqn:Mf})  are row equivalent.
Both ILPs~(\ref{ilp:BF}) and~(\ref{eqn:Mf}) have the same set of non-negative variables without having any other inequality constraints. Hence, the LP relaxations of both ILPs  have the same feasible set and objective function,   and consequently the same symmetry group $G(k,2,t)^{\rm LP}$. Now, if $t$ is even, then by Lemma~\ref{lem:contains} the formulation symmetry group of ILP~(\ref{eqn:Mf}) contains $G(k)^{\rm OD}$.
Since the formulation symmetry group of ILP~(\ref{eqn:Mf}) is a subgroup of   $G(k,2,t)^{\rm LP}$, $G(k,2,t)^{\rm LP}$ contains $G(k)^{\rm OD}$.   The converse statement follows from Lemma~\ref{lem:containsnot}.\qed
\end{pf} 
 Now, inequality~(\ref{eqn:Arquette}) follows from Lemmas~\ref{lem:wr} and~\ref{lem:contains2}. Theorem~\ref{thm:use} 
 follows from Remark~\ref{rem:sol}, Lemmas~\ref{lem:wr} and~\ref{lem:contains2}, and comparing group sizes.
By taking $\x=(N/s^k)\1$, we see that the LP relaxation of ILP~(\ref{ilp:BF}) is feasible, so Method \ref{meth:Glp} of Section~\ref{sec:orbitproj} applies. 
\makeatletter
\renewcommand*\env@matrix[1][*\c@MaxMatrixCols c]{%
  \hskip -\arraycolsep
  \let\@ifnextchar\new@ifnextchar
  \array{#1}}
\makeatother

The following lemma identifies an LP in standard form that has the same feasible set as the LP relaxation of ILP~(\ref{ilp:BF}).
\begin{lem}\label{lem:std}
Let $1\leq t\leq k-1$, $p_{\rm max}\neq \lambda/(s^{k-t})$, and $\A'(k,s,t)\x=\bb'(k,s,t)$ be the equality constraints  of   ILP~(\ref{ilp:BFimproved}), where $p_{\rm max}\leq \lambda$ is the upper bound for
 the variables in ILP~(\ref{ilp:BF}).
 Let
  
\[
\begin{bmatrix}[r|r]
 \B'(k,s,t) & \ddd'(k,s,t)
\end{bmatrix} =
\left\{
\begin{array}{ll}
       \multirow{2}{*}
{$
\begin{bmatrix}[r|r]
-\I & \0  
\end{bmatrix}
$}
  & \text{if $x_1\geq 0$ is a facet and $x_1\leq p_{\rm max}$ is not a facet of}\\
  &\text{the  LP relaxation of  ILP~(\ref{ilp:BF})}, \\
 \multirow{2}{*} {$\begin{bmatrix}[r|r] 
\I & p_{\rm max}\1  
\end{bmatrix}$} & \text{if $x_1\leq p_{\rm max}$ is a facet and $x_1\geq 0$ is not a facet of} \\  
&\text{the LP relaxation of
ILP~(\ref{ilp:BF})}, \\[1ex]
  \begin{bmatrix}[r|r] 
-\I \ &  \0\\
\I \ &   p_{\rm max}\1
\end{bmatrix} & {otherwise. }
\end{array}\right.\]
  Then 
 \begin{equation}\label{eqn:either}
\begin{array}{rl}
&\quad \quad \quad \ \min \   \1^{\top}\x \\
\mbox{s.t.} &  \A'(k,s,t) \x=\bb'(k,s,t), \\
&  \B'(k,s,t)\x\leq \ddd'(k,s,t) 
\end{array}
\end{equation}
is an LP in standard form with the same feasible set and objective function as the LP relaxation of ILP~(\ref{ilp:BF}). 
 \end{lem}
 \begin{pf}
Let $\A(k,s,t)\x=\bb(k,s,t)$ be the equality constraints of ILP~(\ref{ilp:BF}). It is easy to see that the equality constraints $\A(k,s,t)\x=\bb(k,s,t)$ can be obtained as linear combinations of the equality constraints  $\A'(k,s,t)\x=\bb'(k,s,t)$ and vice versa.  
By  Lemma 2 in~\cite{Rosenberg1995},  $${\rm rank}(\A(k,s,t))=\sum_{i=0}^{t} {k \choose i}=m,$$ and $m$ is also equal to the number of rows of $\A'(k,s,t)$. Hence 
 $\A'(k,s,t)$ has full row rank. 
Then it suffices to show that either each inequality in $\0\leq\x$, and/or each inequality in $\x\leq p_{\max}\1$  is a facet of the LP relaxation of ILP~(\ref{ilp:BF}). 

Let LP~(\ref{ilp:BF}) be the LP relaxation of ILP~(\ref{ilp:BF}) and $\mathcal{F}$ be its feasible set.
Since $G(k,s,t)$ preserves $\mathcal{F}$,
 it sends facets of LP~(\ref{ilp:BF}) to its facets and preserves the set of all equality constraints 
satisfied by each point in $\mathcal{F}$. Moreover,  $G(k,s,t)$ acts transitively 
on the variables $x_i$. Thus, it acts transitively on the inequalities $0\leq x_i$ as well as on $x_i\leq p_{\max}$.
So, if one of the  inequalities is satisfied as an equality by each point in $\mathcal{F}$, then either $x_i=0$ for $i \in \{1,\ldots,s^k\}$ or
 $x_i=p_{\max}$ for $i \in \{1,\ldots,s^k\}$ for each  $\x \in \mathcal{F}$. We cannot have $x_i=0$ for $i \in \{1,\ldots,s^k\}$ or $x_i=p_{\max}$ for $i \in \{1,\ldots,s^k\}$ when $p_{\max} \neq \lambda/(s^{k-t})$ as such points are not in $\mathcal{F}$. 
 Therefore, no inequality of LP~(\ref{ilp:BF}) can be satisfied as an equality by each point in $\mathcal{F}$. 
Moreover, at least one of the inequalities of LP~(\ref{ilp:BF}) is a facet. 
Otherwise the feasible set would not be bounded. Then, since $G(k,s,t)$ acts transitively on the variables $x_i$ and preserves the set of facets of LP~(\ref{ilp:BF}),
  $0\leq x_i$ for $i \in \{1,\ldots,s^k\}$ and/or $x_i\leq p_{\max}$ for $i \in \{1,\ldots,s^k\}$ are all facets of LP~(\ref{ilp:BF}). 
\qed
 \end{pf}
Let $S_{s^k}$ be the group of all permutations of coordinates of vectors in $\mathbb{R}^{s^k}$. By Lemma~\ref{lem:std} and the fact that  $(N/s^k)\1\in \mathcal{T}^{\rm LP(\ref{eqn:either})}_{{\rm Fix}_{H}}$  for any subgroup $H$ of $S_{s^k}$, we get 
$$G(k,s,t)^{\rm LP} =
G^{\rm Null}_{(\A'(k,s,t),\B'(k,s,t), \ddd'(k,s,t),\1)}= G_{\PP_{\A'(k,s,t)^{\top}}},
$$ and Theorem~\ref{thm:Glp} follows. 

Method~\ref{meth:GPAT} can be used to compute the automorphism group $G_{\PP_{\A'(k,s,t)^{\top}}}$ of $\PP_{\A'(k,s,t)^{\top}}$ by taking $\A:=\A'(k,s,t)$, $\cc:=\1$ as inputs and stopping once the Step~\ref{step:edgecolored} computation finishes. Then, the output is $H_{\PP_{\A^{\top}}}=G_{\PP_{\A^{\top}}}=G_{\PP_{\A'(k,s,t)^{\top}}}$.
By using Method~\ref{meth:GPAT},  we computed $G(k,s,t)^{\rm LP}=G_{\PP_{\A'(k,s,t)^{\top}}}$  for many $k,s,t$ combinations.
Our computational results and Theorems~\ref{thm:Arquette1} 
and~\ref{thm:Arquette} suggest that  inequality~(\ref{eqn:Arquette}) and the containments in 
Lemmas~\ref{lem:wr} 
 and~\ref{lem:contains2}  are in fact equalities for $ 1 \leq t \leq k-1$ unless $s=2$ and $k=t+1$.   
(In~\cite{Arquette2016}, it was  proven that $G(3,2,2)^{\rm LP}\cong (S_4 \times S_4)\rtimes S_2 $ 
by using {\tt GAP} and Method~\ref{meth:GPAT} of Section~\ref{sec:orbitproj}, where $|G(3,2,2)^{\rm LP}|=1,152>2^34!=192$. By using 
Method~\ref{meth:GPAT} of Section~\ref{sec:orbitproj}, we also observed that $|G(t+1,2,t)^{\rm LP}|>2^{t+1}(t+2)!$ for $t=3,\ldots,10$ and $(|G(t+1,2,t)^{\rm LP}|)/(2^{t+1}(t+2)!)$ increases exponentially with $t$.) In fact, the only known cases for which $1 \leq t \leq k-1$ and yet inequality~(\ref{eqn:Arquette}) is not satisfied as an equality are when $k=t+1$ and $s=2$.

We excluded the case $t=k$ from our results. This is because this case is trivial and completely solved by the following remark.
\begin{rem}
When $t=k$, $x_i=\lambda$ for $i \in \{1,\ldots,s^k\}$ is the unique solution to ILP~(\ref{ilp:BF}). Consequently, the symmetry group of 
ILP~(\ref{ilp:BF}) in this case is $S_{s^k}$, where
$S_{s^k}$ is the set of all permutations of the variables (frequencies) in ILP~(\ref{ilp:BF}).  
\end{rem}
\section{Computational experiments}\label{sec:comparisons}
A speed comparison of exploiting $G(k,s,t)$ and $G(k,s,t)^{\rm LP}=G_{\PP_{\A'(k,s,t)^{\top}}}$ for ILP~(\ref{ilp:BFimproved}) and $G(k,s,t)^{\rm LP\leq}$ for the Theorem 1 ILP in~\cite{Geyer2014}  within  B\&B with isomorphism pruning~\cite{Margot2007} is made in Table~\ref{tbl:FormComp}. The groups $G(k,s,t)$ and 
$G(k,s,t)^{{\rm LP}\leq}$ were computed by using the method in~\cite{Margot2010,Rehn2017}
 as formulation symmetry groups of ILP~(\ref{ilp:BF}) and Theorem 1 ILP in~\cite{Geyer2014}.
 The group $G(k,s,t)^{\rm LP}=G_{\PP_{\A'(k,s,t)^{\top}}}=H_{\PP_{\A'(k,s,t)^{\top}}}$ was computed by using Method~\ref{meth:GPAT} as described at the end of Section~\ref{sec:OAapply}. The automorphism group in each of these methods was computed by using 
 {\tt nauty 25.1}~\cite{McKay2013,McKayP}. 
A computer program written in {\tt C} was used for computing $nce$ in Step~\ref{step:nce} and constructing the  edge colored graph between Step~\ref{step:edgestart} and Step~\ref{step:edgend} in Method~\ref{meth:GPAT}.  
 A singular value decomposition $\U\D\V^{\top}$ and
  $\PP_{\A'(k,s,t)^{\top}}=\V\I_n^{(p)}\V^{\top}$ 
  in Step~\ref{step:svd} 
    were computed in {\tt MATLAB 8.0}~\cite{Matlab:2012}. 
{\tt ISOP 1.1} implementation~\cite{Margot2007} that calls the {\tt CPLEX 12.5.1} libraries~\cite{Cplex} was used for B\&B with isomorphism pruning. The overall running times and the numbers of non-isomorphic solutions pertaining to exploiting $G(k,s,t)$ and $G(k,s,t)^{\rm LP\leq}$  except the OA$(160,8,2,4)$ and OA$(176,8,2,4)$  cases (second, fourth, fifth, and seventh columns) are copied from~\cite{Geyer2014}.
 All cases were run on an \texttt{HP Z820} workstation with
64GB of RAM and a 3.10 GHz \texttt{Intel(R) Xeon(R) E5-2687W} processor.  (Processor information in exploiting $G(k,s,t)$ and $G(k,s,t)^{\rm LP\leq}$ for the results in~\cite{Geyer2014} that we provide here was not provided in~\cite{Geyer2014}.)
 For each $\oan$, the second and the third columns report the number of non-isomorphic solutions enumerated for  ILP~(\ref{ilp:BFimproved}) using $G(k,s,t)$  and $G(k,s,t)^{\rm LP}$. (These are also the number of non-isomorphic solutions of ILP~(\ref{ilp:BF}) using $G(k,s,t)$ and $G(k,s,t)^{\rm LP}$.) The fourth column reports the number of non-isomorphic  solutions
found from the Theorem~1 ILP formulation  in~\cite{Geyer2014} using $G(k,s,t)^{\rm LP\leq}$. The fifth, sixth, and the seventh columns report the times it took to find all 
non-isomorphic solutions using $G(k,s,t)$, 
$G(k,s,t)^{\rm LP}$, and $G(k,s,t)^{\rm LP\leq}$ with  ILP~(\ref{ilp:BFimproved}), ILP~(\ref{ilp:BFimproved}), and the
Theorem~1 ILP formulation  in~\cite{Geyer2014}. Each of these times includes the time it took to compute the exploited symmetry group. The times in parentheses, on the other 
hand, are the times needed to compute the corresponding symmetry groups. 
For most $N,k,s,t$ cases in Table~\ref{tbl:FormComp}, the time needed to find all non-isomorphic solutions  is much greater than that for 
computing the corresponding symmetry groups. 

{ \centering
\tiny
\begin{longtable}{ccccccc}
\caption{Speed comparisons and the number of non-isomorphic solutions}\label{tbl:FormComp}\\
\hline
&{ILP~(\ref{ilp:BFimproved})}&{ILP~(\ref{ilp:BFimproved}) }&
{ILP in~\cite{Geyer2014}}&{ILP~(\ref{ilp:BFimproved}) }&{ILP~(\ref{ilp:BFimproved})  }&{ILP in~\cite{Geyer2014}}\\
{OA($N,k,s,t$)}&{$G(k,s,t)$ }&{$G(k,s,t)^{\rm LP }$ }&{$G(k,s,t)^{\rm LP \leq}$}&{ $G(k,s,t)$}& {$G(k,s,t)^{\rm LP }$ }& {$G(k,s,t)^{\rm LP \leq}$} \\
&{ \# of OAs}& { \# of OAs}&{\# of OAs}& {Times (sec.)}& { Times (sec.)}
&{ Times (sec.)}\\
\hline
\endfirsthead

\hline \multicolumn{5}{c}%
{\tablename\ \thetable{} -- {\rm continued \ from \ previous \ page}}\\
\hline
&{ILP~(\ref{ilp:BFimproved})}&{ILP~(\ref{ilp:BFimproved})}&
{ILP in~\cite{Geyer2014}}&{ILP~(\ref{ilp:BFimproved}) }&{ILP~(\ref{ilp:BFimproved})  }&{ILP in~\cite{Geyer2014}}\\
{OA($N,k,s,t$)}&{$G(k,s,t)$ }&{ $G(k,s,t)^{\rm LP }$}&{$G(k,s,t)^{\rm LP \leq}$}&{$G(k,s,t)$ }& {$G(k,s,t)^{\rm LP }$ }& {$G(k,s,t)^{\rm LP \leq}$} \\
&{ \# of OAs}& { \# of OAs}&{\# of OAs}& { Times (sec.)}& {Times (sec.)}
&{ Times (sec.)}\\
\hline
\endhead
\hline \multicolumn{5}{r}{{\rm Continued \ on \ next \ page}}\\
\hline
\endfoot
 \hline
\endlastfoot
\hline
OA(20,6,2,2) & 75 & 23 & 3,069 & 1 (0) & 7 (6) & 64 (6)\\
OA(20,7,2,2) & 474 & 102 & 51,695 & 13 (0) & 9 (6)& 2,578 (7)\\
OA(20,8,2,2) & 1,603 & 211 & 383,729 & 109 (1)& 22 (7) & 66,377 (11)\\
OA(20,9,2,2) & 2,477 & 351 & 1,157,955 & 485 (4)& 67 (14)& 879,382 (26)\\
OA(20,10,2,2) & 2,389 & 260 & $\geq28,195$ & 1,684 (33)& 215 (72)& $\geq 37,214 $ (76)\\ \hline
OA(24,5,2,2) & 63 & 31 & 723 & 1 (0) & 10 (6)& 18 (6)\\
OA(24,6,2,2) & 1,350 & 274 & 62,043 & 22 (0) & 12 (6)& 1,381 (6)\\
OA(24,7,2,2) & 57,389 & 7,990 & 6,894,001 & 1,721 (0)& 257 (6) & 428,220 (7)\\
OA(24,8,2,2) & 1,470,157 & 165,596 & 4,505,018 & 99,738 (1) & 10,082 (7) & 653,671 (11)\\
OA(24,9,2,2) & 3,815,882 & 1,309,475 & - & 763,643 (4)& 223,138 (14)& - (25)\\
OA(24,5,2,3) & 1 & 1 & 2 & 0 (0)& 6 (6)& 12 (6)\\
OA(24,6,2,3) & 2 & 2 & 5 & 0 (0)& 7 (6)& 12 (6)\\
OA(24,7,2,3) & 1 & 1 & 5 & 0 (0)& 9 (6)& 16 (8)\\
OA(24,8,2,3) & 1 & 1 & 6 & 1 (1)& 14 (7)& 23 (13)\\
OA(24,9,2,3) & 1 & 1 & 6 & 6 (4)& 26 (13)& 44 (30)\\
OA(24,10,2,3) & 1 & 1 & 5 & 55 (42)& 104 (57)& 129 (91)\\
OA(24,11,2,3) & 1 & 1 & 3 & 520 (359)& 540 (441)& 461 (320)\\ 
\hline
OA(32,6,2,3) & 10 & 10 & 31 & 2 (0)& 8 (6)& 12 (6)\\
OA(32,7,2,3) & 17 & 17 & 76 & 2 (0)& 8 (6)& 16 (8)\\
OA(32,8,2,3) & 33 & 33 & 194 & 7 (1) & 14 (7)& 77 (13)\\
OA(32,9,2,3) & 34 & 34 & 364 & 24 (5)& 33 (13)& 658 (30)\\
OA(32,10,2,3) & 32 & 32 & 561 & 102 (42)& 112 (56)& 7,338 (91)\\
OA(32,11,2,3) & 22 & 22 & $\geq 441$ & 560 (364)& 597 (442) & $\geq 36,463$ (319)\\ \hline
OA(40,6,2,3) & 9 & 9 & 65 & 1 (0) & 7 (6) & 13 (6)\\
OA(40,7,2,3) & 25 & 25 & 580 & 2 (0) & 9 (6)& 41 (8)\\
OA(40,8,2,3) & 105 & 105 & 6,943 & 20 (1)& 27 (7)& 4,178 (13)\\
OA(40,9,2,3) & 213 & 213 & 43,713 & 206 (5)& 215 (13)& 260,919 (30)\\
OA(40,10,2,3) & 353 & 353 & $\geq 1,511$ & 1,765 (42)& 1,694 (57)& $\geq 36,279$ (91)\\ \hline
OA(48,7,2,3) & 397 & 397 & 13,469 & 34 (0) & 40 (7) & 862 (8)\\
OA(48,8,2,3) & 8,383 & 8,383 & 896,963 & 2,232 (1)& 2,237 (8) & 552,154 (13)\\
\hline
OA(54,5,3,3) & 4 & 4 & 49 & 2 (1)& 10 (7)& 36 (13)\\
OA(54,6,3,3) & 0 & 0 & 0 & 17 (13)& 37 (24)& 167 (53)\\
 \hline
OA(56,6,2,3) & 86 & 86 & 1,393 & 4 (0)& 11 (6) & 36 (6)\\
OA(56,7,2,3) & 4,049 & 4,049 & 285,184 & 443 (0)& 450 (6)& 20,415 (8)\\ \hline
OA(64,7,2,4) & 7 & 4 & 21 & 99 (0)& 260 (6)& 15 (8)\\
OA(64,8,2,4) & 3 & 2 & 10 & 12 (1)& 38 (8)& 23 (14)\\ \hline
OA(80,6,2,4) & 1 & 1 & 6 & 1 (0)& 7 (6)& 12 (7)\\
OA(80,7,2,4) & 0 & 0 & 0 & 0 (0)& 8 (7)& 15 (8)\\ \hline
OA(81,5,3,4) & 1 & 1 & 2 & 16 (1)& 23 (7)& 20 (13)\\ \hline
OA(96,7,2,4) & 4 & 2 & 31 & 3 (0)& 10 (6)& 15 (8)\\
OA(96,8,2,4) & 0 & 0 & 0 & 2 (1)& 11 (8)& 60 (15)\\ \hline
OA(112,6,2,4) & 3 & 2 & 25 & 1 (0)& 8 (6)& 13 (6)\\ 
OA(112,7,2,4) & 0 & 0 & 0 & 1 (0)& 8 (6)& 18 (8)\\\hline
OA(144,8,2,4) & 20 & 7 & 3,392 & 1,793 (1)& 774 (8)& 1,535,314
 (14) \\ \hline
OA(160,8,2,4) & 99,618 & 11,712 & - & 123,180 (1)& 32,880 (9)& -
(14) \\ \hline
OA(176,8,2,4) & 1,157,443 & 129,138 & - &- (1)& 1,067,822 (8)& - (14)\\ \hline
OA(162,6,3,4) & 0 & 0 & 0 & 20 (14)& 32 (24) & 267 (62)
\end{longtable}
}

The set of all non-isomorphic solutions under the action of
$G(k,s,t)$ and $G(k,s,t)^{\rm LP}$ correspond to a set of all non-isomorphic and non-OD-equivalent $\oan$.  The numbers of all non-isomorphic solutions obtained by exploiting $G(k,s,t)$ and $G(k,s,t)^{\rm LP}$   for the bottleneck cases 
OA$(160,8,2,4)$ and OA$(176,8,2,4)$ corroborate the numbers of all 
non-isomorphic and non-OD-equivalent OA$(160,8,2,4)$ and OA$(176,8,2,4)$ in~\cite{Bulutoglu2016}.
The number of  all non-isomorphic solutions under the action of $G(k,s,t)^{\rm LP\leq}\cong S_{s-1}\wr S_k$ equals  the number of all $\oan$ up to a weaker form of isomorphism. 

 For cases in which $G(k,s,t)^{\rm LP}$ captures symmetries not in $G(k,s,t)$, the speedup gleaned from adding slack variables to the Theorem~1 ILP formulation in~\cite{Geyer2014} and using $G(k,s,t)^{\rm LP}$ to enumerate $\oan$ up to OD-equivalence under the action of $G(k,s,t)^{\rm LP}$ with ILP~(\ref{ilp:BFimproved}) becomes the fastest enumeration method for OAs as the number of variables $s^k$ increases. (ILP~(\ref{ilp:BFimproved}) can be obtained from the Theorem 1 ILP formulation in~\cite{Geyer2014} by adding slack variables and the $x_1 \geq 1 $ inequality.) Moreover, this 
 speedup appears to grow exponentially with the number of variables. However, the cases OA$(64,7,2,4)$  and 
OA$(24,11,2,3)$ are exceptions to this trend. Hence, exploiting a larger symmetry group  drastically overcomes the extra computational burden due to having additional variables. This underscores the importance of developing tools for finding larger subgroups of the symmetry group of an ILP. Finally, 
 the cost of computing $G(k,2,2)^{\rm LP}$ when $k\geq 6$ and $G(k,2,4)^{\rm LP}$ when $k \geq 8$ is more than compensated for with the speedup gleaned from exploiting the additional symmetries not in $G(k,2,2)$ and $G(k,2,4)$.
For many $s=2$ and even $t$ cases and all the bottleneck cases, using the larger symmetry group $G(k,2,t)^{\rm LP}$ drastically reduces solution times.

A set of all non-isomorphic (non-OD-equivalent if $s=2$ and $t$ is even) $\oan$ can be obtained by adding columns to a set of all non-isomorphic (non-OD-equivalent if $s=2$ and $t$ is even) OA$(N,k-1,s,t)$. Bulutoglu and Ryan~\cite{Bulutoglu2016} used this fact to develop the Hybrid method that enumerates a set of all non-isomorphic (non-OD-equivalent) $\oan$ by adding columns to a set of all non-isomorphic (non-OD-equivalent) OA$(N,k-1,s,t)$.  This method adds columns to input 
OA$(N,k-1,s,t)$ by finding a set of all non-isomorphic solutions  to ILPs derived from 
the input OA$(N,k-1,s,t)$ and ILP~(\ref{ilp:BFimproved}). For each input OA$(N,k-1,s,t)$, it uses B\&B with isomorphism pruning with  a group depending on the input OA$(N,k-1,s,t)$.
 However, this method removes only some of the symmetry within the B\&B with isomorphism pruning algorithm and requires converting $\oan$ to graphs and using {\tt nauty}~\cite{McKay2013,McKayP}  for removing isomorphic graphs that correspond to isomorphic (OD-equivalent) $\oan$~\cite{Bulutoglu2016,Ryan2010}.
 
McKay~\cite{McKay1998,McKay1996} had previously developed a technique for generating combinatorial objects with partial isomorph rejection 
when it is possible to sequentially obtain larger objects from the smaller. In particular, McKay's technique is applicable to generating a set of all non-isomorphic (non-OD-equivalent) $\oan$ from a set of all non-isomorphic (non-OD-equivalent) OA$(N,k-1,s,t)$. However, 
 just like the Hybrid method,  when this technique is applied to the problem of generating a set of all non-isomorphic (non-OD-equivalent) $\oan$ it does not completely eliminate the need to use {\tt nauty} for removing isomorphic (OD-equivalent) $\oan$. 
 In fact, Bulutoglu and Ryan~\cite{Bulutoglu2016} implemented McKay's technique for generating a set of all non-OD-equivalent OA$(160,8,2,4)$ and OA$(176,8,2,4)$ from a set of all non-OD-equivalent OA$(160,7,2,4)$ and OA$(176,7,2,4)$ 
 and observed that the running times for the Hybrid method were $1/(11.44)$ and $1/(1.44)$ times those of McKay's technique. The OA$(160,8,2,4)$ and OA$(176,8,2,4)$ are the largest $2$-symbol, strength $4$ OAs that have been classified, where the use of a symmetry exploiting method was necessary~\cite{Bulutoglu2016}.

Unlike the Hybrid method or McKay's technique, exploiting $G(k,2,t)^{\rm LP}$ when $t$ is even within B\&B with isomorphism pruning  enabled us to directly generate a set of all non-OD-equivalent
$\oant$ without using  {\tt nauty}~\cite{McKay2013,McKayP} to remove OD-equivalent OAs. However, we did use {\tt nauty}~\cite{McKay2013,McKayP}  to find $G(k,s,t)^{\rm LP}$.
This was a viable method by Theorem~\ref{thm:use}, and it reduced the enumeration times  of a set of all non-OD-equivalent OA$(160,8,2,4)$ and OA$(176,8,2,4)$ in comparison to the Hybrid method in~\cite{Bulutoglu2016}  by  factors of $1/(2.16)$ and $1/(1.36)$. (We ran the OA$(176,8,2,4)$ case 
on our \texttt{HP Z820} workstation with
$128$GB of RAM and $2.00$ GHz \texttt{Intel(R) Xeon(R) E5-2650} processor as well to allow making  comparisons to the corresponding times for the Hybrid method and McKay's technique in~\cite{Bulutoglu2016}.) 
Hence, using $G(k,2,t)^{\rm LP}$ as described in this paper reduced the running time for finding all OD-equivalence classes of OA$(160,8,2,4)$ and OA$(176,8,2,4)$ by  factors of $1/(24.71)$ and $1/(1.96)$ in comparison to McKay's technique.
\section{Conclusion }\label{sec:future}
In this paper, we showed that there may be hidden symmetries in an LP that cannot be captured by the
formulation symmetry group. These symmetries are either masked by redundant constraints or due to equality
constraints. As a remedy, we developed a method that captures all the symmetries of a feasible LP. (The symmetry group of an infeasible LP is isomorphic to $S_n$.) We tested our method on the LP relaxations of a family of ILPs for classifying  OAs, and for $\oant$ with even $t$, we found LP relaxation symmetry groups with drastically larger sizes than their corresponding formulation symmetry groups.  Finally, we exploited the newly found larger groups $G(k,2,t)^{\rm LP}$ within B\&B with isomorphism pruning. This enabled us to improve the times it took to find all OD-equivalence classes of OA$(160,8,2,4)$ and OA$(176,8,2,4)$ by  factors of $1/(2.16)$ and $1/(1.36)$.

 One of the key findings of this article involves the enumeration of a set of
 all non-isomorphic solutions to an ILP. In this context,  converting the inequality constraints to equalities by introducing slack variables and  using the LP relaxation $\Glp$ of the resulting ILP within B\&B with isomorphism pruning can reduce the enumeration  time by several orders of magnitude. In particular, this method would be most useful in determining whether a given ILP is feasible. We propose testing this idea along with the methods in this paper on the MIPLIP problems studied in~\cite{Liberti2012,Rehn2017} as a future research project. A limited preliminary study on the MIPLIB problems
 in~\cite{Liberti2012}  suggests that, when computing $\Glp$, generating the graphs between Step~\ref{step:edgestart} and Step~\ref{step:edgend} in Method~\ref{meth:GPAT} and computing the double coset decompositions in Method~\ref{meth:Glpp} are the major bottlenecks in terms of both time and memory requirements. Based on our experience with the OA problem, we expect that the  time requirements will be much greater 
 for finding  sets of all non-isomorphic optimal solutions of the MIPLIB problems under the action of their respective LP relaxation symmetry groups than that for computing their LP relaxation symmetry groups.
\section*{Acknowledgements}
The authors thank two anonymous referees for improving the paper. The authors also thank Mr. David Doak for general computer support.
This research was supported by the AFOSR grant F4FGA04013J001. 
The views expressed in this article are those of the authors and do not reflect the official policy or position of the United States Air Force, Department of Defense, or the U.S. Government.\\
\bibliographystyle{elsarticle-harv}
\bibliography{bibliography}
\end{document}